\documentclass[a4paper,12pt]{amsart}
\usepackage{amssymb}
\usepackage{amsmath}
\def\R{\mathbb{R}}
\def\C{\mathbb{C}}

\def\so{\mathfrak{so}}
\def\sp{\mathfrak{sp}}

\def\m{\mathfrak{m}}
\def\SO{\mathrm{SO}}
\def\Spin{\mathrm{Spin}}
\def\SU{\mathrm{SU}}
\def\Sp{\mathrm{Sp}}
\def\T{\mathrm{T}}
\def\TM{\mathrm{T}M}

\def\Cas{\mathrm{Cas}\,}
\def\End{\mathrm{End}\,}
\def\Hom{\mathrm{Hom}\,}

\def\id{\mathrm{id}}

\def\Ric{\mathrm{Ric}\,}
\def\scal{\mathrm{scal}\,}
\def\Ad{\mathrm{Ad}\,}
\def\#{\sharp}

\def\G{\mathrm{G}}

\def\Aut{\mathrm{Aut}}

\def\bea{\begin{eqnarray*}}
\def\eea{\end{eqnarray*}}
\def\Ric{\mathrm{Ric}}
\def\tr{\mathrm{tr}}
\def\la{\langle}
\def\ra{\rangle}

\def\GL{\mathrm{GL}}
\def\Spin{\mathrm{Spin}}
\def\SO{\mathrm{SO}}
\def\O{\mathrm{O}}
\def\Id{\mathrm{Id}}
\def\<#1,#2>{\langle\,#1,\,#2\,\rangle}
\newtheorem{Lemma}{Lemma}[section]
\newtheorem{Proposition}[Lemma]{Proposition}
\newtheorem{Theorem}[Lemma]{Theorem}
\newtheorem{Corollary}[Lemma]{Corollary}
\newtheorem{Definition}[Lemma]{Definition}
\newtheorem{Remark}[Lemma]{Remark}
\def\proof{\noindent\textbf{Proof:}\quad}

\def\qed{\quad\hfill\ensuremath{\Box}}

\def\be{\begin{equation}}
\def\ee{\end{equation}}
\newcommand{\pref}[1]{(\ref{#1})}

\begin{document}
%
%
\title{Deformations of nearly parallel $\G_2$-structures}
\author{B. Alexandrov, U. Semmelmann} 
\thanks{The first author was partially supported by Contract 195/2010 with Sofia University "St. Kl. Ohridski".}
\address{Bogdan Alexandrov\\
Faculty of Mathematics and Informatics\\
University of Sofia 'St. Kliment Ohridski'\\
James Bourchier Blvd. 5 \\
1164 Sofia, Bulgaria}
\email{balexand@fmi.uni-sofia.bg}
\address{Uwe Semmelmann\\
Institut f\"ur Geometrie und Topologie \\
Fachbereich Mathematik\\
Universit{\"a}t Stuttgart\\
Pfaffenwaldring 57 \\
70569 Stuttgart, Germany
}
\email{uwe.semmelmann@mathematik.uni-stuttgart.de}
\begin{abstract}
We study the infinitesimal deformations of a proper nearly parallel
$\G_2$-structure and prove that they are characterized by a certain
first order differential equation. In particular we show that the
space of infinitesimal deformations modulo the group of
diffeomorphisms is isomorphic to a subspace of co-closed
$\Lambda^3_{27}$-eigenforms of the Laplace operator for the
eigenvalue $8\scal/21$. We give a similar description for the space
of infinitesimal Einstein deformations of a fixed nearly parallel
$\G_2$-structure. Moreover we show that there are no deformations on
the squashed $S^7$ and on $\SO(5)/\SO(3)$, but that there are
infinitesimal deformations on the Aloff-Wallach manifold
$N(1,1)=\SU(3)/U(1)$.
\end{abstract}
\maketitle
\section{Introduction}
A nearly parallel $\G_2$-structure on a $7$-dimensional manifold $M$ is
given by a $3$-form  $\sigma$ of special algebraic type satisfying the
differential equation $\ast  d\sigma = \tau_0 \sigma$ for some constant $\tau_0$. Such a manifold has
a structure group contained in the exceptional Lie group $\G_2 \subset
\SO(7)$ and, in particular, a Riemannian metric $g$ induced by
$\sigma$. It can be shown that nearly parallel $\G_2$-manifolds are
irreducible and Einstein with scalar curvature $\scal = \tfrac{21}{8}\,\tau^2_0$. Moreover, the existence
of such a structure is equivalent to the existence of a spin structure
with a Killing spinor.

\medskip

Another equivalent description of nearly parallel $\G_2$-structures is in
terms of the metric cone $(\hat M, \hat g)$, which has to have
holonomy contained in $\Spin(7)$, considered as subgroup of $\SO(8)$.
The metric cone is the manifold $\hat M = M \times \R_+$ with the
warped product metric $\hat g =  r^2 g \oplus dr^2$. If $(M^7,g)$ is simply connected and
not isometric to the standard sphere, then there are three possible
cases: the holonomy of $(\hat M, \hat g)$ is contained in $Sp(2)$,
equivalently, $(M^7,g)$ is a $3$-Sasakian manifold, the holonomy can be
$\SU(4)$, equivalently, $(M^7,g)$ is an Einstein-Sasaki manifold, or
the holonomy is precisely $\Spin(7)$, in which case we call the
$\G_2$-structure {\it proper}. We recall that these three cases
correspond to the existence of a 3-, 2- resp. 1-dimensional space of
Killing spinors. Proper nearly parallel $\G_2$-structures are also
characterized by the vanishing of the Lie derivative $L_\xi \sigma$
for any Killing vector field $\xi$.

\medskip

In this article we shall mainly consider the case of proper nearly
parallel $\G_2$-manifolds. In \cite{uwe1} it is shown that any
$7$-dimensional $3$-Sasakian manifold admits a second nearly parallel
$\G_2$-structure which is proper. The corresponding Einstein metric belongs to the
metrics of the canonical variation of the $3$-Sasakian Einstein
metric. Applying this construction to the homogeneous $3$-Sasakian spaces $S^7$ and
$N(1,1)$ one obtains homogeneous proper nearly parallel $\G_2$-structures: the squashed
$7$-sphere and the second Einstein metric on $N(1,1)$. The
Aloff-Wallach spaces $N(k,l)$ for $(k,l)\neq (1,1)$ also have exactly two nearly parallel
$\G_2$-structures, both of which are proper.  A further example is the isotropy
irreducible space $\SO(5)/\SO(3)$. In fact, due to the
classification~\cite{uwe1} these are the only
homogeneous nearly parallel $\G_2$-manifolds.

\medskip

As a last remarkable property of nearly parallel $\G_2$-manifolds we
mention the existence of a metric connection $\bar \nabla$ with
totally skew-symmetric torsion. The so-called {\it canonical
connection} $\bar\nabla$ is defined as $\bar\nabla = \nabla -
\tfrac{\tau_0}{12}\sigma$ and has holonomy contained in the group $\G_2
\subset \SO(7)$. Nearly parallel $\G_2$-manifolds appear as one of
two exceptional cases in a classification of metric connections with
parallel torsion due to Cleyton and
Swann~\cite{cleyton-swann}. The other exceptional case is the class
of 6-dimensional nearly K\"ahler manifolds, which turns out to be in
various ways rather similar to nearly parallel $\G_2$-manifolds. The
defining condition is the existence of a nearly parallel almost
complex structure $J$, i.e., $J$ satisfying $(\nabla_X J) (X)=0$ for
any vector field $X$. Nearly K\"ahler manifolds in dimension 6 are
also Einstein manifolds admitting a Killing spinor. Moreover, the
metric cone has holonomy contained in $\G_2$.

\medskip

In this article we shall show that nearly parallel $\G_2$-manifolds
are also in another respect very similar to nearly K\"ahler
manifolds: the description of infinitesimal deformations.
In~\cite{uwe2} the space of infinitesimal nearly K\"ahler deformations is
identified with the space of primitive co-closed $(1,1)$-eigenforms
of the Laplace operator for the eigenvalue $2\scal/5$,
\cite{andrei-uwe} contains a similar description of the space of
infinitesimal Einstein deformations. This space turns out to be the
sum of three such eigenspaces. Finally, in \cite{andrei-uwe2} it is
shown that infinitesimal deformations for the known homogenous
examples only exist in the case of the flag manifold $\SU(3)/\T^2$.
For all three results we shall 
obtain a
counterpart on nearly parallel $\G_2$-manifolds.

\medskip

We start with the equations of R.~Bryant (cf.
Proposition~\ref{bryant} and \cite{bryant}) describing the
infinitesimal deformation of an arbitrary $\G_2$-structure. They give
equations for the tangent vector on a curve of $\G_2$-structures.
Specializing to the case of nearly parallel
$\G_2$-structures and staying transversal to the action of the
diffeomorphism group, we obtain that the space of such deformations
is a direct sum of two spaces, $D_1$ and $D_3$, consisting of $1$-forms
and $3$-forms respectively. As shown in Section~\ref{secD1}, the space
$D_1$ parametrizes Einstein-Sasakian structures compatible with the given nearly parallel
$G_2$-structure. The more interesting space is $D_3$ which consists of the solutions
$\phi$ in $\Lambda^3_{27}\T^*M$ of the
differential equation $\ast d \phi = -\tau_0 \phi$. In particular, infinitesimal
deformations $\phi \in D_3$ are co-closed and eigenforms of the Hodge-Laplace operator
for the eigenvalue $\tau^2_0=8\scal/21$. But more important for the computation in
examples is that they are also eigenforms for the eigenvalue $\frac{5\tau_0^2}{6}$
of the $\G_2$-Laplace operator $\bar\Delta$ introduced in Section~\ref{laplace}. In
Section~\ref{sec5} we describe the space of infinitesimal Einstein
deformations of the metric of a nearly parallel
$G_2$-structure. In addition to $D_3$ one obtains two other spaces
of sections of $\Lambda^3_{27}\T^*M$ which are characterized by similar
equations. In the last section we compute the infinitesimal Einstein
deformations of the normal homogeneous examples: the isotropy
irreducible space $\SO(5)/\SO(3)$, the squashed $7$-sphere and the second Einstein metric
on the Aloff-Wallach space $N(1,1)$. We show that there exist no Einstein deformations and,
in particular, no deformations of the nearly parallel $G_2$-structure in the first two cases,
while in the third the space of infinitesimal Einstein deformations coincides with the space
of infinitesimal nearly parallel $G_2$-deformations and is $8$-dimensional. We do not know
whether these infinitesimal deformations integrate to real Einstein deformations.

%
\section{Preliminaries}
Let $e_1,\ldots, e_7$ denote the standard basis of $\R^7$ and
$e^1,\ldots,e^7$ its dual
basis. On $\R^7$ we fix the canonical scalar
product $\la \cdot, \cdot \ra$  and the standard orientation. We
shall write $e^{i_1\ldots i_k}$ for the wedge product $e^{i_1}\wedge
\ldots \wedge e^{i_k}\in \Lambda^k (\R^7)^*$ and define the {\it
fundamental 3-form} as
\begin{equation}\label{sigma}
\sigma = e^{123} + e^{145} + e^{246} + e^{347} - e^{167} + e^{257} -
e^{356}.
\end{equation}
The exceptional group $\G_2$ is defined as the subgroup of
$\GL(7,\R)$ that fixes the 3-form $\sigma$, i.e., $\G_2 = \{ g \in
\GL(7,\R) \,\, | \,\, g^*\sigma = \sigma\}$. The group $\G_2$ is a
14-dimensional compact, connected, simple Lie group, which acts
irreducibly on $\T:=\R^7$ and preserves the metric, the orientation and the
Hodge dual of $\sigma$, i.e. the 4-form
\begin{equation}\label{*sigma}
\ast \sigma = e^{4567} +  e^{2367} +  e^{1357} +  e^{1256} -  e^{2345}
+
e^{1346} -  e^{1247}.
\end{equation}

The irreducible representations of $\G_2$ can be indexed by their
highest weights, which are pairs of non-negative integers $(p,q)$ if written as
linear combinations of the two fundamental weights. The
corresponding representation will be denoted by $V_{p,q}$. In this
paper we will in particular be interested  in the following four irreducible
$\G_2$-representations: the trivial representation $V_{0,0} =
\mathbb{R}$, the standard representation $V_{1,0} =
\T:=\mathbb{R}^7$, the adjoint representation $V_{0,1} =
\mathfrak{g}_2$ and the representation on traceless symmetric 2-forms $V_{2,0}=
S_0^2 \T^*$. Among the irreducible representations
these are uniquely determined by their dimensions $1$, $7$, $14$ and
$27$ respectively. Therefore we shall use the dimensions as lower
indices when we decompose the space of $k$-forms $\Lambda^k \T^*$ into
irreducible components. In other words, $\Lambda^k_r$ will denote the
$r$-dimensional irreducible subspace of $\Lambda^k \T^*$. With this
notation we have
\begin{equation}\label{split}
\Lambda^2 = \Lambda^2\T^* = \Lambda^2_7 \oplus \Lambda^2_{14}, \qquad
\Lambda^3 = \Lambda^3\T^* = \Lambda^3_1 \oplus \Lambda^3_7 \oplus \Lambda^3_{27},
\end{equation}
with an isomorphic decomposition for $\Lambda^4\T^* \cong \Lambda^3\T^*$
and
$\Lambda^5\T^* \cong
\Lambda^2\T^*$ obtained with the help of the Hodge $*$-operator.
The one-dimensional spaces in $\Lambda^3$ resp. $\Lambda^4$ are spanned by $\sigma$ resp.
$\ast \sigma$. The space $\Lambda^2_{14}$ is isomorphic to the Lie algebra of $\G_2$ and
the other subspaces can be characterized by
$$
\Lambda^2_7 = \{ X \lrcorner \, \sigma \in \Lambda^2 \,\, | \,\, X
\in \T\} \cong \T, \quad \Lambda^3_7 = \{ X \lrcorner \ast\sigma \in \Lambda^3
\,\, | \,\, X \in \T\} \cong \T,$$
$$\Lambda^3_{27} = \{\alpha \in \Lambda^3
\,\, | \,\, \alpha \wedge \sigma = 0 = \alpha \wedge \ast \sigma\} \cong V_{2,0}.
$$
 In the sequel we shall use the following $\G_2$-equivariant
 isomorphisms, which were introduced by Bryant in \cite{bryant}: ${\bf i}: S_0^2 \T^* \to \Lambda^3_{27}$ and ${\bf j}:
 \Lambda^3_{27} \to S_0^2\T^*$, where ${\bf i}$ is the restriction to
$S_0^2 \T^* \subset S^2 \T^*$ of the map $S^2 \T^* \to \Lambda^3\T^*$,
 defined on decomposable elements by
$$
\alpha \odot \beta \mapsto \alpha \wedge (\beta \lrcorner \sigma) +
\beta \wedge (\alpha \lrcorner \sigma),$$
while ${\bf j}$ is given by
$${\bf j}(\gamma)(X,Y) = *((X \lrcorner \sigma) \wedge (Y \lrcorner \sigma) \wedge \gamma).$$
Note that ${\bf j} = -8{\bf i}^{-1}$. With the help of ${\bf i}$ one
can obtain
explicit elements of $\Lambda^3_{27}$, e.g.
\begin{equation}\label{expl}
{\bf i}(e^1 \odot e^2) = e^{146} + e^{157} + e^{245} - e^{267}.
\end{equation}
\noindent
Because of $\T^* \otimes \T^* = S^2 \T^* \oplus \Lambda^2 \T^*$
we have the following decomposition:
\begin{equation}\label{dec1}
V_{1,0} \otimes V_{1,0} \cong \mathbb{R} \oplus V_{2,0} \oplus V_{1,0} \oplus V_{0,1}.
\end{equation}
Later we shall also need the decompositions
\begin{equation}\label{dec2}
V_{1,0} \otimes V_{2,0} \cong V_{1,0} \oplus V_{2,0} \oplus V_{0,1} \oplus V_{1,1} \oplus V_{3,0},
\end{equation}
\begin{equation}\label{dec3}
V_{1,0} \otimes V_{0,1} \cong V_{1,0} \oplus V_{2,0} \oplus V_{1,1}.
\end{equation}

The group $\G_2$ can also
be defined as the stabilizer of the vector cross product $P$, given by
\begin{equation}\label{sigma1}
\sigma(X,Y,Z) = \la P(X,Y), Z \ra ,
\end{equation}
where $X,Y,Z$ are any vectors in $\T$. Recall from \cite{gray-fernandez} that a $2$-fold vector cross product $P$ is a bilinear map $P: \T \times \T \rightarrow \T$ satisfying for
all $X,Y \in \T $ the equations
\begin{equation}\label{P1}
\la P(X,Y),X \ra = \la P(X,Y),Y\ra = 0 \quad \mbox{and}\quad \|P(X,Y)\|^2 =  \|X \|^2 \|Y \|^2 -\la X,Y \ra^2.
\end{equation}
In particular, it follows from the second equation of~(\ref{P1}) that
$P$ is skew-symmetric. Thus we can consider $P$ as a linear map $P :
\Lambda^2\T \rightarrow \T$ and write $P(X\wedge Y) = P(X,Y)$. In this
notation the second equation of~(\ref{P1}) reads: $\|P(X\wedge
Y)\|^2=\|X\wedge Y\|^2$. We also refer to \cite{gray-fernandez} for
the following relations satisfied by a general 2-fold vector cross product:
\begin{Lemma}\label{P2}
For $X, Y, Z \in \T$ we have
$$
\begin{array}{ll}
& (1) \quad   \la P(X,Y), Z\ra= \la X,P(Y,Z)\ra  ,\\
& (2) \quad   P(X,P(X,Y))  = - \|X\|^2Y + \la X,Y\ra X  ,\\
& (3) \quad   2 P(P(X,Y),Z) =
P(P(Y,Z),X) + P(P(Z,X),Y)   
+ 3 \la X,Z\ra Y -  3 \la Y,Z\ra X  .
\end{array}
$$
\end{Lemma}

\noindent
From now on we will usually identify vectors and 1-forms via the
metric and denote with $\{e_i\}$, $i=1,\ldots,7$ an orthonormal basis
of $\T$. For later use we still note
\begin{Lemma} Let $X$ and $Y$ be any vectors in $\T$. Then the
  following equations hold
\begin{eqnarray}
&
(X \,\lrcorner\, \sigma) \wedge \sigma \;=\; -2\, X \wedge \,\ast\,\sigma, \label{ai1}
\\
&
(X \,\lrcorner\, \sigma) \wedge \ast \sigma \;=\; 3\,\ast X  , \label{ai2}
\\
&
\sum_i\,(e_i \,\lrcorner\,X \,\lrcorner\,\sigma)\,\lrcorner\, (e_i \wedge \sigma)
\;=\; 3\,X\,\lrcorner\,\ast\sigma  , \label{ai3}
\\
&
\sum_i\,(e_i \,\lrcorner\,X \,\lrcorner\,\sigma) \wedge (e_i \lrcorner\,\sigma)
\;=\; 3\,X\,\lrcorner\,\ast\sigma  , \label{ai4}
\\
&
(X\,\lrcorner\,Y \,\lrcorner\,\sigma) \,\lrcorner\, \sigma \;+\; X\,\lrcorner\,Y \,\lrcorner\, * \sigma
\;=\; - X\wedge Y, \label{ai5}
\\
&
P(X \,\lrcorner\,\sigma) \;=\; 3 X. \label{ai6}
\end{eqnarray}
\end{Lemma}

The $\GL_7$-orbit of $\sigma$ in $\Lambda^3\T^*$ is an open
set by dimensional reasons. As usual it is denoted with $\Lambda^3_+$. Forms in $\Lambda^3_+$
are called {\it stable} or {\it definite}.

Let $M$ be a 7-dimensional manifold. The union of the subspaces $\Lambda^3_+ \T_x^*M$, $x \in M$, of stable
forms defines an open subbundle $\Lambda^3_+ \T^*M \subset \Lambda^3\T^*M$. There is a one-to-one
correspondence between $\G_2$-structures on $M$, i.e.  reductions of the structure group of $M$ to
the group $\G_2$, and the space of sections of $\Lambda^3_+ \T^*M$, which we will
denote with $\Omega^3_+(M)$. The defining 3-form $\sigma \in \Omega^3_+(M)$
determines a Riemannian metric $g$ and an orientation of $M$ via the relation
\begin{equation}\label{metric}
-6g(X,Y)*1 = X \lrcorner \sigma \wedge Y \lrcorner \sigma \wedge \sigma
\end{equation}
($*1$ denotes the volume form). Let $\nabla$ be the Levi-Civita
connection of $g$. Then the covariant derivative $\nabla \sigma$ is a
section of the bundle $\T^*  \otimes \mathfrak{g}_2^\perp $, where $\mathfrak{g}_2^\perp
\cong \T^*$ is the orthogonal complement of $\mathfrak{g}_2$ in
$\Lambda^2 \T^*$ and we identify bundles with the $\G_2$-representation
defining it. It follows from  \pref{dec1} that this bundle decomposes
as $\mathbb{R} \oplus V_{2,0} \oplus V_{1,0} \oplus V_{0,1}$ and thus
the covariant derivative of $\sigma$ has four components.
Accordingly, one has the 16 Fernandez-Gray classes of $\G_2$-structures,
with the four basic classes $W_1, W_2, W_3,  W_4$ corresponding to the
four irreducible summands.

In this article we shall consider the class $W_1$ of so called
{\it nearly parallel} (or {\it weak})
$\G_2$-structures, i.e. $\G_2$-structures induced by a non-parallel
3-form $\sigma \in \Omega^3_+(M)$, such that $\nabla \sigma$
is a section of the 1-dimensional subbundle defined by the trivial
$\G_2$-representation. Nearly parallel $\G_2$-structures can be
described by several equivalent conditions in terms of $\sigma$.

\begin{Proposition}\label{np}
Let $M$ be a 7-dimensional manifold with a $\G_2$-structure defined by
a
3-form $\sigma \in \Omega^3_+(M)$. Then the following
conditions are equivalent
\begin{enumerate}
\item The 3-form $\sigma$ defines a nearly parallel $\G_2$-structure.
\item The 3-form $\sigma$ is a Killing 3-form, i.e. $\nabla \sigma =
\tfrac14 d\sigma $.
\item There exists a $\tau_0 \in \R \setminus \{ 0 \}$ with $\nabla
\sigma = \tfrac{\tau_0}{4} \ast \sigma $.
\item There exists a $\tau_0 \in \R \setminus \{ 0 \}$ with $\nabla_X
(\ast \sigma) = - \tfrac{\tau_0}{4} X \wedge  \sigma$ for all
vector fields $X$.
\item There exists a $\tau_0 \in \R \setminus \{ 0 \}$ with $d \sigma
  =\tau_0
\ast \sigma$.
\item $X \lrcorner \nabla_X \sigma = 0$ holds for all vector fields $X$.
\end{enumerate}
\end{Proposition}
\proof The equivalence of (3) and (4) is obvious, while the
equivalence of
(1), (2), (3) and (6) has been proved in \cite{gray-fernandez}. The
only
point not mentioned there is that $\tau_0$ is constant. This fact is
also
known (see e.g. \cite{uwe1}) and can be proven as follows. Since (5)
is an
obvious consequence of (3), we can differentiate it to obtain $d\tau_0
\wedge *\sigma = 0$, which implies $d\tau_0 = 0$. Finally, that (5)
implies the remaining conditions was proved in \cite{uwe1}. This
is the only point where one uses that $\tau_0$ is different from
zero.
\qed

Let $P$ be the associated vector cross product, defined in
(\ref{sigma1}).
Then the condition (6) of the proposition above is equivalent to
$(\nabla_X P)(X, Y) = 0$ for any vector fields $X,Y$, i.e., to $P$
being {\it nearly parallel} \cite{gray}. Further straightforward
consequences of Proposition~\ref{np} in the case of nearly parallel
$\G_2$-manifolds are: $d^* \sigma = 0$ and
$\Delta \sigma = \tau^2_0 \,\sigma$, where here and in the
following $\Delta=dd^* + d^* d$ denotes the Hodge-de Rham Laplacian.
Moreover it follows that $\sigma$ is a special Killing 3-form,
i.e. the additional equation
$\nabla_X d\sigma = - \frac14 \tau_0^2 X^* \wedge \sigma $ is satisfied
for all vector
fields $X$ (cf. \cite{uwe}).

The canonical connection $\bar\nabla$ of a $\G_2$-structure is the
unique $\G_2$-connection whose torsion is equal to the intrinsic torsion of
the $\G_2$-structure.
In  the nearly parallel case it has totally skew-symmetric and
parallel torsion and is explicitly given by
\begin{equation}\label{nablabar}
g(\bar\nabla_XY, \,Z) \;=\; g(\nabla_XY,\,Z)
\;-\;\tfrac{\tau_0}{12}\,\sigma(X,\,Y,\,Z)
\end{equation}
or, equivalently, by
\begin{equation}\label{nablabar2}
\bar\nabla_X = \nabla_X - \tfrac{\tau_0}{12} P_X,
\end{equation}
where the endomorphism $P_X$ is defined by $P_X Y := P(X,Y)$.

\begin{Remark}\label{r1}
{\rm The fact that $P$ is $\G_2$-invariant allows the following important
application, which we shall use several times in this article. Let $V$
be an irreducible $\G_2$-representation contained in some tensor space and
$VM$ be the corresponding associated bundle. Then the endomorphism
$P_X$ extends to an endomorphism of $VM$ and we may consider the
$\G_2$-equivariant map $V \to \T^* \otimes V$, defined by
$
\varphi \mapsto \sum_i e^i \otimes P_{e_i} \varphi,
$
which we again denote by $P$. By \pref{nablabar2} we have
$$(\bar\nabla - \nabla)\varphi = - \tfrac{\tau_0}{12} P\varphi$$
for any section $\varphi$ of $VM$. Let $U$ be an irreducible
component
in $\T^* \otimes V$. Suppose first that $U$ is not isomorphic to
$V$ as a $\G_2$-representation. Then there exists no non-zero
$\G_2$-equivariant map from $V$ to $U$ and therefore the
$UM$-part $(P\varphi)_{UM}$ of $P\varphi$ vanishes, which implies
$(\nabla \varphi)_{UM} = (\bar\nabla \varphi)_{UM}$. On the other
hand, if
$U$ is isomorphic to $V$, then $U = i(V)$, where $i:V \to \T^* \otimes
V$
is some $\G_2$-equivariant embedding. Let $\pi:\T^* \otimes V \to U$ be
the
projection. Then $\pi \circ P: V \to U$ is also $\G_2$-equivariant and
therefore by Schur's lemma $\pi \circ P = ci$ for some constant $c$.
Thus $(\nabla \varphi)_{i(VM)} = (\bar\nabla \varphi)_{i(VM)} +
\tfrac{c\tau_0}{12} i(\varphi)$. Finally, since $\bar\nabla \varphi$
and $P\varphi$ are sections of $\T^*M \otimes VM$, the same is true
for $\nabla \varphi$, despite the fact that $\nabla$ is not a
$\G_2$-connection.}
\end{Remark}

\begin{Remark}\label{r2}
{\rm Our choice of the orientation induced by a stable $3$-form $\sigma$ is the opposite of the choice of Bryant in \cite{bryant}. As a consequence our $*$, ${\bf j}$, $\tau_0$ and $f_1$ from the next section differ from those in \cite{bryant} by a sign. }
\end{Remark}

%
\section{Deformations of $\G_2$-structures}\label{sec3}
In this section we will consider a smooth curve $\sigma_t$ of
nearly parallel $\G_2$-structures and describe its tangent vector
$\dot \sigma$ in $t=0$. Here and in the sequel the dot denotes the
time derivative at $t=0$. As a starting point we use the
following result of R.~Bryant \cite{bryant}
for curves of arbitrary $\G_2$-structures (cf. also \cite{Joyce}).

\begin{Proposition}\label{bryant}
Let $(M^7,g)$ be a Riemannian manifold with a family $\sigma_t \in
\Omega^3_+(M)$ of $\G_2$-structures. Let $g_t$ be the family of
metrics and $\ast_t$ the Hodge star operator associated with $\sigma_t$. Then
there exist three time-dependent differential forms
$f_0 \in \Omega^0(M)$, $f_1 \in \Omega^1(M)$ and $f_3 \in
\Omega^3_{27}(M)$ that satisfy the
equations
\begin{enumerate}
\item \;
$\dot \sigma \;=\; 3\,f_0\,\sigma \;+\; \ast(f_1 \wedge \sigma)
\;+\; f_3  ,$
\item \;
$\dot g \;=\; 2\,f_0\,g \;-\; \tfrac12\,{\bf j}(f_3)  ,$
\item \;
$\dot{\ast \sigma} \;=\; 4\,f_0\,\ast\sigma \;+\; f_1 \wedge \sigma
\;-\; \ast  f_3  ,$
\item \;
$\dot{\ast 1} \;=\; 7\,f_0\,\ast 1 .$
\end{enumerate}
\end{Proposition}

Our aim is to study deformations of a given nearly parallel
$\G_2$-structure $\sigma$ on a compact manifold $M$ by nearly parallel
$\G_2$-structures $\sigma_t$. We will only be interested in
deformations of the nearly parallel $\G_2$-structures modulo the
action of the group $\mathbb{R}^* \times
\mathrm{Diff}(M)$, given by $$(\lambda,\varphi) \cdot \sigma = \lambda^3
\varphi
(\sigma) = \lambda^3 (\varphi^{-1})^* \sigma \circ \varphi.$$
If $\sigma$ induces the metric $g$, the Hodge dual $*\sigma$ and
the volume form $*1$, then $\widetilde{\sigma} = \lambda^3 \sigma$
induces
\begin{equation}\label{scale}
\widetilde{g} = \lambda^2 g, \quad \widetilde{*}\widetilde{\sigma} =
\lambda^4 *\sigma, \quad \widetilde{*}1 = \lambda^7 *1, \quad
\widetilde{\tau}_0 = \frac{1}{\lambda}\tau_0.
\end{equation}
Therefore we can always assume that the volume of $M$ with respect
to $g$ is normalized. Moreover, we can apply the Ebin's Slice
Theorem and assume that $g_t$ is a curve in the slice through
$g$. A nearly parallel $\G_2$-structure is Einstein with scalar curvature
\begin{equation}\label{scal}
\scal_{g} \;=\; \tfrac{21}{8}\,\tau_0^2 .
\end{equation}
Thus $\dot{g}$ is an infinitesimal Einstein deformation of $g$ and by
the
theorem of Berger-Ebin (see \cite{besse}, Chapter 12) we have
\begin{equation}\label{ide}
\tr \, \dot g \;=\; 0, \qquad\delta \,\dot g \;=\;0, \qquad \Delta_L
\dot g = \tfrac{2\scal}{7} \dot g \;=\; \tfrac{3}{4}\,\tau_0^2 \dot g,
\end{equation}
where $\Delta_L$ is the Lichnerowicz Laplacian (see \cite{besse} or
Section~\ref{sec5} below).
Since $\tr \,\dot g = 14\,f_0$, it immediately follows that $f_0$
vanishes and the equations of Proposition~\ref{bryant} may be
rewritten as
\begin{equation}\label{deform}
\dot \sigma \;=\;\ast(f_1 \wedge \sigma) \;+\; f_3,\quad \dot g\;=\;
-\tfrac12\,{\bf j}(f_3),\quad
 \dot{\ast \sigma} \;=\; f_1\wedge \sigma \;-\; \ast  f_3,\quad
  \dot{\ast1} \;=\; \,0,
\end{equation}
while equations \pref{ide} become
\begin{equation}\label{ide0}
\delta \,{\bf j}(f_3) \;=\;0, \qquad \Delta_L \, {\bf j}(f_3) =
\tfrac{3}{4}\,
\tau_0^2 \, {\bf j}(f_3).
\end{equation}
The fact that $\sigma_t$ is a family of nearly parallel
$\G_2$-structures
means by definition that
\begin{equation}\label{npt}
d\sigma_t = \tau_0(t)\,\ast \sigma_t
\end{equation}
for some function $\tau_0(t)$. However, $g_t$ is a family of Einstein
metrics and therefore $\scal_{g_t}$ is constant as function in $t$, as follows from
Corollary~2.12 of \cite{besse}. This, together with
\pref{scal}, implies that the function $\tau_0$ is constant too. Thus,
differentiating \pref{npt} with respect to $t$, we obtain the
linearized equation $d\dot \sigma = \tau_0 \, \dot{\ast \sigma},$
which by \pref{deform} yields
\be\label{first}
d\ast(f_1\wedge \sigma) \;+\;df_3 \;=\; \tau_0 (f_1 \wedge \sigma
\;-\; \ast f_3)  .
\ee

\noindent
The discussion above motivates the following definition.

\begin{Definition}
{\rm An {\it infinitesimal (nearly parallel) deformation} of a compact
nearly parallel $\G_2$-manifold $(M,\sigma)$ is a section $(f_1,f_3)$
of the bundle $\Lambda^1\T^*M \oplus \Lambda^3_{27}\T^*M$, which
satisfies the equations from \pref{ide0} and \pref{first}.
}
\end{Definition}
The rest of this section is devoted to deriving a more explicit
description of the space of infinitesimal deformations of nearly
parallel $\G_2$-structures. In a first step we obtain information
about the 3-form component $f_3$ of the infinitesimal deformation.
\begin{Lemma}\label{lemf3}
The covariant derivatives of $f_3$ with respect to $\bar \nabla$ and
$\nabla$ have no component in $\T\subset \T^* \otimes
\Lambda^3_{27}$, i.e., it holds  $(\bar\nabla f_3)_\T = (\nabla f_3)_\T
= 0$.
In particular, the differential and codifferential of $f_3$ satisfy
the equations:
$$
(df_3)_{\Lambda^4_{1}} = 0, \qquad
(df_3)_{\Lambda^4_{7}} = 0,\qquad (d^*f_3)_{\Lambda^2_{7}}=0 ,\qquad
(d*f_3)_{\Lambda^5_7}=0 .
$$
\end{Lemma}
\proof The divergence $\delta$
is defined as the composition of the covariant derivative $\nabla$
and the equivariant contraction $c:\T^* \otimes S^2_0\T^* \rightarrow
\T^*$. It follows from Remark~\ref{r1} and the decomposition
\pref{dec2}
that $(\nabla {\bf j}(f_3))_\T = (\bar \nabla {\bf j}(f_3))_\T$ and
$(\nabla f_3)_\T = (\bar \nabla f_3)_\T$. This implies
\begin{equation}\label{deltaj}
\delta {\bf j}(f_3) = -c \nabla {\bf j}(f_3) = -c (\nabla
{\bf j}(f_3))_\T = -c (\bar \nabla {\bf j}(f_3))_\T =
-c \circ (1\otimes {\bf
  j})
(\bar\nabla f_3)_\T ,
\end{equation}
where the last equality follows from ${\bf j}$ being $\G_2$-equivariant
and $\bar\nabla$-parallel. Since $c \circ (1\otimes {\bf j})$ is non-zero on
$\T\subset \T^* \otimes
\Lambda^3_{27}$ (which can be checked on an explicit element), we
finally obtain $\delta {\bf j}(f_3) = 0$ if and only if $(\nabla
{\bf j}(f_3))_\T = (\bar\nabla f_3)_\T = 0$.

By the definition of the differential we have $df_3 =
\varepsilon \nabla f_3$, with the $\G_2$-equivariant wedging map
$\varepsilon : \T^* \otimes \Lambda^3 \T^* \rightarrow \Lambda^4
\T^*$. Again by \pref{dec2} and Remark~\ref{r1} we see that $\nabla f_3$ has
no components in bundles associated with the trivial representation $\mathbb{R}$.
Thus $df_3$ has no components in $\Lambda^4_{1}$, which proves the first
equation for $df_3$. The second follows from $(df_3)_{\Lambda^4_{7}} = \varepsilon (\nabla f_3)_\T = 0$ and the
remaining two are proved in a similar way.
\qed

\noindent
In the next step we will derive information about the 1-form part
of infinitesimal deformations.

\begin{Proposition}\label{c2}
For the 1-form $f_1$ the following holds:
\begin{enumerate}
\item $\bar \nabla f_1 = -\frac{1}{3} \,\tau_0\, f_1 \,\lrcorner\, \sigma $.
\item $\nabla f_1 = -\frac{1}{4} \,\tau_0\, f_1 \,\lrcorner\, \sigma $.
\item $f_1$ is a Killing $1$-form and
$df_1= -\tfrac12 \,\tau_0\, f_1 \,\lrcorner\, \sigma $. In
particular, $d^*f_1 = 0$ and $(df_1)_{\Lambda^2_{14}}=0$.
\item
$d(f_1 \wedge \ast \sigma) = -\tfrac{3}{2} \tau_0 \ast f_1$, \quad in
particular \quad $d^*(f_1\,\lrcorner\,\sigma )  = -\tfrac32 \tau_0
f_1$.
\item $\Delta f_1 =
\tfrac34 \tau_0^2 f_1 $.
\item
$d^* (f_1\wedge\sigma )  =  -\tau_0 f_1 \,\lrcorner \, \ast
\sigma$. In particular \quad $d\ast (f_1\wedge \sigma) =
\tau_0 f_1 \wedge \sigma$.
\item $f_1$ has constant length.
\end{enumerate}
\end{Proposition}
\proof Statements (2) and (3) are obviously equivalent and (1) and (2) are
equivalent
because \pref{nablabar} implies
$
\bar \nabla f_1 = \nabla f_1 -\frac{\tau_0}{12} f_1 \lrcorner \sigma
$.
The remaining properties are
consequences of each of the first three.
To prove (6) we use $d^* f_1 =
0$,
$d^* \sigma = 0$, (3) of Proposition~\ref{np} and \pref{ai4} to obtain
$$
d^*(f_1 \wedge \sigma) = \sum \nabla_{e_i}
f_1 \wedge (e_i \lrcorner \sigma ) - \nabla_{f_1^\sharp} \sigma
= -\frac{1}{4} \tau_0 \sum (e_i \lrcorner f_1 \lrcorner \sigma) \wedge
(e_i \lrcorner \sigma ) - \frac{1}{4} \tau_0 f_1 \lrcorner *\sigma =
-\tau_0 f_1 \lrcorner *\sigma  .
$$
Property (7) follows from
$$
d|f_1|^2 (X) = 2 \langle \nabla_X f_1, f_1 \rangle =
-\frac{1}{2}\tau_0
\langle X \lrcorner f_1 \lrcorner \sigma, f_1 \rangle =
-\frac{1}{2}\tau_0
f_1 \lrcorner X \lrcorner f_1 \lrcorner \sigma = 0.
$$

For the proof of (3) we start with computing the
$\Lambda^2_7$ part of $df_1$. Since $\tau_0$ is a non-zero constant
we may take the differential of equation~(\ref{first})
to obtain $d(f_1\wedge \sigma)=d\ast f_3$. Now, let
$\beta$ be the $1$-form, defined by $(df_1)_{\Lambda^2_7} = \beta
\lrcorner \,\sigma $. Then the vanishing of the $\Lambda^5_7$-part
of  $d*f_3$ and equation~\pref{ai1} imply
$$
0 = (d(f_1\wedge \sigma))_{\Lambda^5_7} = (df_1)_{\Lambda^2_7} \wedge
\sigma - \tau_0 f_1 \wedge \ast \sigma
= (\beta \lrcorner \sigma) \wedge \sigma - \tau_0 f_1 \wedge \ast
\sigma = (-2\beta - \tau_0 f_1) \wedge \ast \sigma  .
$$
Since the wedge
product with $\ast \sigma$ defines an injective map on 1-forms, we
find $\beta = -\tfrac12 \tau_0 f_1$. Hence $(df_1)_{\Lambda^2_7} =
-\tfrac12 \tau_0 f_1 \lrcorner \sigma$.
The $\Lambda^2_7$-part of $\nabla f_1$ with respect to the
decomposition \pref{dec1} is $\frac{1}{2}
(df_1)_{\Lambda^2_7}$ and we obtain $(\nabla f_1)_{\Lambda^2_7} = -\tfrac14 \tau_0 f_1 \lrcorner \sigma$.

We continue with proving (4). Since the fundamental 3-form $\sigma$ is
coclosed, it follows $d \ast \sigma =0$. Moreover, the wedge product
with $\ast \sigma$ defines an equivariant map $\Lambda^2 \rightarrow
\Lambda^6\cong\Lambda^1$, which by Schur's Lemma vanishes on
$\Lambda^2_{14}$. Hence, using \pref{ai2} we obtain
$$
d(f_1\wedge \ast \sigma) = df_1\wedge \ast \sigma = (df_1)_{\Lambda^2_7}\wedge \ast \sigma = -\tfrac12 \tau_0 (f_1\lrcorner \sigma)\wedge \ast \sigma
= -\tfrac32 \tau_0 \ast f_1.
$$
As an
immediate consequence we obtain in addition
$d^*(f_1\,\lrcorner\,\sigma ) = -\tfrac32 \tau_0 f_1$. Since $\tau_0 \neq 0$, this implies $d^*f_1 =
0$.

Next we want to show that the $\Lambda^3_7$-part of $d( f_1 \lrcorner \sigma)$ vanishes. Using Proposition~\ref{np} we compute
$$
\begin{array}{rl}
d( f_1 \lrcorner \sigma) &= \sum_i e^i \wedge \nabla_{e_i}( f_1 \lrcorner \sigma)
\;=\; \sum_i e^i \wedge (\nabla_{e_i} f_1 \lrcorner \sigma + f_1 \lrcorner \nabla_{e_i}\sigma)\\[1ex]
&= -\sum_i \nabla_{e_i} f_1  \lrcorner (e_i \wedge \sigma) +(d^* f_1)\sigma - f_1 \lrcorner d\sigma + \nabla_{f_1^\sharp}\sigma \\[1ex]
&= -\Phi(\nabla f_1) -
\tfrac34 \tau_0 f_1 \lrcorner *\sigma,
\end{array}
$$
where $\Phi: \Lambda^1 \otimes \Lambda^1 \to \Lambda^3$ denotes the map
$\Phi(\gamma) = \sum_i (e_i \lrcorner \gamma) \lrcorner (e^i \wedge \sigma)$.
Since $\Phi$ is obviously $G_2$-invariant, we obtain
$$(d( f_1 \lrcorner \sigma))_{\Lambda^3_7} = -(\Phi(\nabla f_1))_{\Lambda^3_7} -
\tfrac34 \tau_0 f_1 \lrcorner *\sigma= -\Phi((\nabla f_1)_{\Lambda^2_7}) -
\tfrac34 \tau_0 f_1 \lrcorner *\sigma = \frac{1}{4}\tau_0 \Phi(f_1 \lrcorner
\sigma) - \tfrac34 \tau_0 f_1 \lrcorner *\sigma.$$


Because \pref{ai3} implies $\Phi(f_1 \lrcorner \sigma) = 3 f_1 \lrcorner *\sigma$, we finally obtain
\begin{equation}\label{f1}
(d( f_1 \lrcorner \sigma))_{\Lambda^3_7} = 0  .
\end{equation}

Now we shall prove the vanishing of the $\Lambda^2_{14}$-part of $df_1$ using the compactness of $M$. From our equation for $(df_1)_{\Lambda^2_7}$ we conclude
$$
0 = d^2 f_1 = d ((df_1)_{\Lambda^2_7} + (df_1)_{\Lambda^2_{14}} ) = -\tfrac12 \tau_0 d(f_1 \lrcorner \sigma) + d (df_1)_{\Lambda^2_{14}}   .
$$
From here it follows with~(\ref{f1}) that $(d (df_1)_{\Lambda^2_{14}} )_{\Lambda^3_7}= 0$. By definition the differential $d$ is the
composition of the invariant wedging map $\varepsilon : T^* \otimes \Lambda^2 \rightarrow \Lambda^3$ and the covariant derivative
$\nabla$. By Remark~\ref{r1} $\nabla \gamma$ is a section of $\T^* M \otimes \Lambda^2_{14} M$ for any section $\gamma$ of $\Lambda^2_{14} M$. Since by \pref{dec3} there is only one component isomorphic to $\T$ in $\T^* \otimes \Lambda^2_{14}$, we obtain for $\gamma := (df_1)_{\Lambda^2_{14}}$:
$$
0 = (d\gamma)_{\Lambda^3_7} = \pi_{\Lambda^3_7} \circ \varepsilon \, \nabla \gamma
= \pi_{\Lambda^3_7} \circ \varepsilon \, (\nabla \gamma)_\T.
$$
Because $\pi_{\Lambda^3_7} \circ \varepsilon$ is different from
zero on $\T \subset \T^*\otimes \Lambda^2_{14}$, as one checks on an explicit element, this yields $(\nabla \gamma)_\T = 0$.

We may use a similar argument for the codifferential $d^*$, which is
the composition of the invariant contraction map $c: \T^* \otimes
\Lambda^2 \rightarrow \Lambda^1$ and the covariant derivative. Hence we have
$$d^*\gamma = - c \nabla \gamma = -c \, (\nabla \gamma)_\T = 0.$$
Then the $L^2$-scalar product of
$d^*\gamma$ and $f_1$ yields
$$
0 = (d^*\gamma, f_1) = (\gamma, df_1) = \|\gamma \|^2  .
$$
Thus it follows that $\gamma = 0$, i.e. $(df_1)_{\Lambda^2_{14}} = 0$,
and
that $df_1$
is indeed a section of $\Lambda^2_7\T^* M$ with $df_1=
(df_1)_{\Lambda^2_7}
 =-\tfrac12 \tau_0
f_1 \,\lrcorner\, \sigma $.

We already know that $f_1$ is coclosed and thus
$$
\Delta f_1 = d^*d f_1 = -\tfrac12 \tau_0 d^*(f_1 \lrcorner
\sigma) =\tfrac34 \tau_0^2 f_1,
$$
which proves (5). Using the fact that the manifold $(M^7, g)$ is Einstein with
$\scal =\tfrac{21}{8}\,\tau_0^2$, we obtain $\Delta f_1 = 2\Ric
(f_1)$.
By the well-known
characterization of Killing vector fields on compact manifolds this
implies that $f_1$ is Killing.
\qed

Finally we combine Proposition~\ref{c2} with the initial equations to
obtain a characterization of infinitesimal deformations of nearly
parallel $\G_2$-structures.

\begin{Theorem}\label{thmid}
The space of infinitesimal deformations of a compact nearly parallel
$\G_2$-manifold $(M,\sigma)$ is the direct sum of the
finite-dimensional spaces
$$
D_1 := \{ f_1 \in \Omega^1(M) \,\, | \,\, \nabla f_1 = -\tfrac{1}{4}
\,\tau_0\, f_1 \,\lrcorner\, \sigma \} \quad \mbox{and} \quad D_3 :=
\{ f_3 \in \Omega^3_{27}(M) \,\, | \,\, *df_3= -\tau_0 f_3 \}.
$$
In particular, $f_1$ and $f_3$ are co-closed eigenforms of the Laplace
operator for the eigenvalues $\frac34\tau^2_0$ and $\tau^2_0$
respectively.
\end{Theorem}
\proof
It remains to prove the equations for $f_3$. For this we substitute
the expression for $d\ast (f_1\wedge \sigma)$ of Proposition~\ref{c2}
back into equation (\ref{first}) and obtain $df_3 = -\tau_0 *f_3$.
Since $\tau_0 \neq 0$ this immediately implies that $f_3$ is coclosed.
Then the Laplace operator is computed as $\Delta f_3 = d^*d f_3 =
(*d)^2 f_3 = \tau_0^2 f_3$. It follows that an infinitesimal
deformation lies in the direct sum of the spaces $D_1$ and $D_3$.
They are finite-dimensional since they are contained in certain
eigenspaces of the Laplace operator.

Conversely, by Proposition~\ref{c2} $\nabla f_1 = -\frac{1}{4}
\,\tau_0\,
f_1 \,\lrcorner\, \sigma$ and $*df_3= -\tau_0 f_3$ imply
\pref{first}.
Further, $df_3= -\tau_0 *f_3$ yields $(df_3)_{\Lambda^4_{7}} = 0$,
i.e., $\pi_{\Lambda^4_7} \circ \varepsilon (\nabla f_3)_\T = 0$,
where $\T\subset \T^* \otimes
\Lambda^3_{27}$ is the component isomorphic to $V_{1,0}$. Since
$\pi_{\Lambda^4_7} \circ \varepsilon$ is non-zero on $\T$ (which
can be checked on an explicit element), it follows that
$(\nabla f_3)_\T = 0$. Thus, by Remark~\ref{r1} also
$(\bar \nabla f_3)_\T = 0$ and therefore by \pref{deltaj}
$\delta {\bf j}(f_3) = -c (1\otimes {\bf j}) (\bar\nabla f_3)_\T = 0$.
It remains to show that $*df_3= -\tau_0 f_3$ implies $\Delta_L {\bf
  j}(f_3)
= \tfrac{3}{4}\,\tau_0^2 {\bf j}(f_3)$, which will be done in
Section~\ref{sec5}. \qed

%
%
\section{$\G_2$-deformations and Sasakian structures} \label{secD1}

In this section we will investigate the relation between nearly parallel $\G_2$-manifolds with a non-trivial space $D_1$ in Theorem~\ref{thmid} and Sasakian structures. The first result in this direction is the following.

\begin{Proposition}\label{propD1}
Let $(M,g,\sigma)$ be a compact nearly parallel $\G_2$-manifold normalized so that $\tau_0 = 4$. Then:
\begin{enumerate}
\item If $\dim D_1 \ge 1$, then $(M,g)$ is a Sasakian-Einstein manifold.
\item If $\dim D_1 \ge 2$, then $(M,g)$ is a 3-Sasakian manifold.
\end{enumerate}
\end{Proposition}
\proof
The assumption about the normalization is not a restriction because of \pref{scale}. Let $0 \neq f_1 \in D_1$, then Proposition~\ref{c2} shows that $f_1$ is a Killing $1$-form of constant length and we can assume $|f_1| = 1$. Thus, to prove that $f_1$ is the contact form of a Sasakian structure it remains (see \cite{BG}) to verify the curvature condition
$$
(\nabla^2_{X,Y} f_1) (Z) \;=\; f_1 (Y) g(X,Z) \,-\, f_1 (Z) g(X,Y)
\;=\; (f_1 \wedge X)(Y,Z)
$$
However, taking the covariant derivative of the defining equation
of $D_1$ immediately implies:
$$
\begin{array}{ll}
(\nabla^2_{X,Y} f_1) (Z)
& = -\frac{1}{4} \,\tau_0\, (\nabla_X (f_1 \,\lrcorner\, \sigma))(Y,Z) \;=\; -(\nabla_X f_1 \,\lrcorner\, \sigma)(Y,Z) \,-\, (f_1  \,\lrcorner\, \nabla_X \sigma)(Y,Z)
\\[1ex]
&= \;\; \,\frac{1}{4} \,\tau_0\, ((X \lrcorner f_1 \,\lrcorner\, \sigma)\,\lrcorner\, \sigma)(Y,Z) \,-\, \frac{1}{4} \,\tau_0 \, (f_1  \,\lrcorner\, X \,\lrcorner\, *\sigma)(Y,Z) \;=\; (f_1 \wedge X)(Y,Z)
\end{array}
$$
where we also used \pref{ai5}. Since $g$ is known to be Einstein, we obtain the first statement.

If $\dim D_1 \ge 2$, then $(M,g)$ has two Sasakian structures, whose contact forms are linearly independent. This implies the second statement (see \cite{BG}, Lemma~8.1.17). \qed

Recall that a $\G_2$-structure on a $7$-dimensional manifold $M$ defines a canonical spin structure on $M$. The $\G_2$-structure is furthermore nearly parallel if and only if the associated spin structure admits real Killing spinors \cite{uwe1}. In this case the nearly parallel $\G_2$-structures inducing the given metric and spin structure are in bijective correspondence with the projectivization of the space of Killing spinors in the real spinor bundle \cite{uwe1}. The complex spinor bundle is the complexification of the real spinor bundle and the space of real Killing spinors is the complexification of the space of Killing spinors in the real spinor bundle, so both spaces have the same dimension over the respective field. After a suitable normalization of the metric (which in our case amounts to ensuring that $\tau_0 =4$) this dimension is also equal to the dimension of the space of parallel spinors on the metric cone $\hat M$ of $M$ for the spin structure induced by the one on $M$. This is a result of B\"ar \cite{Baer} in the simply connected case and holds also in general, as explained by Wang in \cite{Wang}. Hence, as noticed in \cite{Baer}, if $M$ is compact, then either the restricted holonomy group of $\hat M$ is one of $\Spin(7)$, $\SU(4)$, $\Sp(2)$, or $\hat M$ is flat. In the latter case $M$ is a quotient of the standard sphere $S^7$. According to a result of Friedrich \cite{Fried}, all nearly parallel $\G_2$-structures on $S^7$ which induce the standard metric are conjugated under the action of the isometry group. Thus neither $S^7$ nor its quotients admit $\G_2$-deformations. Therefore from now on we shall exclude from our considerations the case of nearly parallel $\G_2$-manifolds with constant curvature. Under this assumption the compact nearly parallel $\G_2$-manifolds split into the following three different types.

{\it Type 1}. The space of real Killing spinors is $1$-dimensional. Then there is only one $3$-form inducing the given metric, orientation and spin structure. We call such nearly parallel $\G_2$-structures {\it proper}. Notice that our definition of a proper nearly parallel $\G_2$-structure is slightly different from those in \cite{uwe1} and \cite{BG}. In \cite{uwe1} one assumes additionally that the manifold is simply connected, while the definition in \cite{BG} requires that the cone has holonomy equal to $\Spin(7)$. For simply connected manifolds the three definitions are equivalent.

{\it Type 2}. The space of real Killing spinors is $2$-dimensional. Then the given metric and orientation are induced by a Sasaki-Einstein structure but not by a $3$-Sasakian structure. In terms of the cone $\hat M$ this is equivalent to saying that the holonomy group of $\hat M$ is contained in $\SU(4)$ but not in $\Sp(2)$. Indeed, the subgroup of $\Spin(8)$ which acts as identity on a $2$-dimensional subspace of one of the half-spin representations is $\Spin(6) \cong \SU(4)$. In this case the $3$-forms inducing the given metric, orientation and spin structure are parametrized by $\mathbb{R}P^1$.

{\it Type 3}. The space of real Killing spinors is $3$-dimensional. The given metric and orientation are induced by a $3$-Sasakian structure. In terms of the cone $\hat M$ this is equivalent to saying that the holonomy group of $\hat M$ is equal to $\Sp(2)$. In this case the $3$-forms inducing the given metric and orientation are parametrized by $\mathbb{R}P^2$.

Now we shall describe the nearly parallel $\G_2$-structures of types 2 and 3 without reference to Killing spinors. Recall that the cone of a Riemannian manifold $(M,g)$ is $(\hat M, \hat g)$, where $\hat M := \mathbb{R}_+ \times M$, $\hat g := dr^2 + r^2g$ and $r$ is the natural coordinate on $\mathbb{R}_+$. As shown in \cite{Baer}, if we normalize the nearly parallel $\G_2$-structure so that $\tau_0 = 4$, then $\sigma = \partial_r \lrcorner \varphi |_{r=1}$, where $\varphi$ is a parallel (and also stable) $4$-form on the cone.

Suppose first that the holonomy group of the cone is equal to $\SU(4)$ (which implies that $M$ is Sasaki-Einstein but not $3$-Sasakian). Then the space of parallel $4$-forms on $\hat M$ is spanned by $\Omega_{\hat I} \wedge \Omega_{\hat I}$, $\mathrm{Re}\, \Psi_{\hat I}$, $\mathrm{Im}\, \Psi_{\hat I}$. Here $\Omega_{\hat I}$ is the K\"ahler form and $\Psi_{\hat I}$ the complex volume form of the $\SU(4)$-structure. Thus
$$\sigma = \partial_r \lrcorner \left. \left( \frac{1}{2}c_0 \Omega_{\hat I} \wedge \Omega_{\hat I} + c_1 \mathrm{Re}\, \Psi_{\hat I} + c_2 \mathrm{Im}\, \Psi_{\hat I} \right) \right|_{r=1}.$$
Equivalently, one can write this as
$$\sigma = c_0 \eta \wedge \Omega + c_1 \mathrm{Re}\, \Psi + c_2 \mathrm{Im}\, \Psi,$$
where $\eta$ is the contact form of the Sasaki-Einstein structure on $M$, $\Omega = \nabla \eta$ is the horizontal K\"ahler form and $\Psi$ is the horizontal complex volume form. Now a straightforward computation using \pref{metric} shows that $\sigma$ induces the given metric and orientation if and only if $c_0 = -1$ and $c_1^2 + c_2^2 = 1$. Hence we have the following explicit $S^1$-family of nearly parallel $\G_2$-structures:
\begin{equation}\label{sigmat}
\sigma_t = -\eta \wedge \Omega + \cos t \, \mathrm{Re}\, \Psi + \sin t \, \mathrm{Im}\, \Psi.
\end{equation}
In particular, each $\sigma_t$ is of type 2.

Now suppose that the holonomy group of the cone is $\Sp(2)$ (i.e., $M$ is a $3$-Sasakian manifold). Then the space of parallel $4$-forms is spanned by
$$\Omega_{\hat I_1} \wedge \Omega_{\hat I_1}, \quad \Omega_{\hat I_2} \wedge \Omega_{\hat I_2}, \quad \Omega_{\hat I_3} \wedge \Omega_{\hat I_3}, \quad \Omega_{\hat I_1} \wedge \Omega_{\hat I_2}, \quad \Omega_{\hat I_2} \wedge \Omega_{\hat I_3}, \quad \Omega_{\hat I_3} \wedge \Omega_{\hat I_1}.$$
Here $\Omega_{\hat I_1}$, $\Omega_{\hat I_2}$, $\Omega_{\hat I_3}$ are the K\"ahler forms of the hyper-K\"ahler structure $\hat I_1$, $\hat I_2$, $\hat I_3$ on the cone (we use the convention $\hat I_1 \hat I_2 = -\hat I_3$). Thus
$$\sigma = \partial_r \lrcorner \left. \left( \frac{1}{2} \sum_\lambda s_{\lambda \lambda} \Omega_{\hat I_\lambda} \wedge \Omega_{\hat I_\lambda} + \sum_{\lambda < \mu} s_{\lambda \mu} \Omega_{\hat I_\lambda} \wedge \Omega_{\hat I_\mu} \right) \right|_{r=1}.$$
Equivalently, this can be written as
\begin{equation}\label{sigmaS}
\sigma = \sum_{\lambda,\mu = 1}^3 s_{\lambda \mu} \eta_\lambda \wedge \Omega_\mu
\end{equation}
with $s_{\lambda \mu} = s_{\mu \lambda}$, where $\eta_1$, $\eta_2$, $\eta_3$ are the contact forms of the $3$-Sasakian structure on $M$ and $\Omega_\lambda = \nabla \eta_\lambda$ are the corresponding K\"ahler forms. Again a straightforward computation using \pref{metric} shows that $\sigma$ given by \pref{sigmaS} induces the given metric and orientation if and only if the matrix $S = \left( s_{\lambda \mu} \right)$ is in $\SO(3)$ and $\mathrm{tr}\, S = -1$. The condition $s_{\lambda \mu} = s_{\mu \lambda}$ means furthermore that $\sigma$ is nearly parallel if and only if $S$ is symmetric. An orthogonal matrix is symmetric if and only if its eigenvalues are real and the condition $\mathrm{tr}\, S = -1$ implies that they are $1,-1,-1$. But an orthogonal matrix with eigenvalues $1,-1,-1$ is completely determined by its $1$-eigenspace. Thus we obtain that the nearly parallel $\G_2$-structures are parametrized by $\mathbb{R}P^2$ (in particular, they are of type 3). We shall identify $\mathbb{R}^3$ with $\mathrm{span}\{ \eta_1, \eta_2, \eta_3 \}$. Then $\eta \in \mathrm{span}\{ \eta_1, \eta_2, \eta_3 \}$ is the contact form of a Saskai-Einstein structure if and only if $\eta$ lies on the unit sphere $S^2$. Let $S(\eta) = \left( s_{\lambda \mu}(\eta) \right)$ denote the orthogonal matrix with eigenvalues $1,-1,-1$ whose $1$-eigenspace is spanned by $\eta$. Then the nearly parallel $\G_2$-structures are
$$\left\{ \sigma_{S(\eta)} = \sum_{\lambda,\mu} s_{\lambda \mu}(\eta) \eta_\lambda \wedge \Omega_\mu \,\, | \,\, \eta \in S^2 \right\}$$
(notice that $S(\eta) = S(-\eta)$).

Fixing an $\eta$, we can again write the $S^1$-family $\sigma_{\eta,t}$ from the $\SU(4)$-case. Inside the $\mathbb{R}P^2$-family it is identified by
$$\{ \sigma_{\eta,t} \} = \{ \sigma_{S(\eta')} \,\, | \,\, \eta' \in S^2, \,\, \eta' \perp \eta \}.$$
This follows from the fact that $\Psi_{\hat I_1} = \frac{1}{2} (\Omega_{\hat I_2} - i\Omega_{\hat I_3}) \wedge (\Omega_{\hat I_2} - i\Omega_{\hat I_3})$, i.e.,
$$\mathrm{Re}\, \Psi_1 = \eta_2 \wedge \Omega_2 - \eta_3 \wedge \Omega_3, \quad \mathrm{Im}\, \Psi_1 = -\eta_2 \wedge \Omega_3 - \eta_3 \wedge \Omega_2.$$

Finally, let the holonomy group $Hol(\hat M)$ of the cone lie strictly between $\SU(4)$ and $\Sp(2)$. Then the restricted holonomy group is $\Sp(2)$ and therefore $Hol(\hat M) \subset \Sp(2)\Sp(1)$ as the normalizer of $\Sp(2)$ in $O(8)$ is $\Sp(2)\Sp(1)$. Now the fact that $Hol(\hat M)$ preserves a complex structure implies $Hol(\hat M) \subset \Sp(2)U(1)$. Finally, $Hol(\hat M)$ preserves a complex volume form, so the $U(1)$ part of $Hol(\hat M)$ is contained in
$$\{ a \in U(1) \,\, | \,\, a^4 = 1 \} = \{ 1,i,-1,-i \} = \mathbb{Z}_4.$$
Since $\Sp(2)\mathbb{Z}_2 = \Sp(2)$, it remains $Hol(\hat M) = \Sp(2)\mathbb{Z}_4$. Now we have to find which $\Sp(2)$-invariant $4$-forms are also $\Sp(2)\mathbb{Z}_4$-invariant. Notice that the action of $i \in \mathbb{Z}_4$ is in fact the complex structure $\hat I = \hat I_1$. Since $\hat I_1$ acts on $\Omega_{\hat I_1}$ as the identity and on $\Omega_{\hat I_2}$ and $\Omega_{\hat I_3}$ as minus identity, the space of $\Sp(2)\mathbb{Z}_4$-invariant $4$-forms is $4$-dimensional and is spanned by
$$\Omega_{\hat I_1} \wedge \Omega_{\hat I_1}, \quad \Omega_{\hat I_2} \wedge \Omega_{\hat I_2}, \quad \Omega_{\hat I_3} \wedge \Omega_{\hat I_3}, \quad \Omega_{\hat I_2} \wedge \Omega_{\hat I_3}.$$
Now the results of the $\Sp(2)$-case imply that $\sigma$ is given by \pref{sigmaS} with $s_{12} = s_{21} = s_{13} = s_{31} = 0$. Thus either $s_{11} = 1$ and
$$\sigma = \sigma_{S(\eta_1)} = \eta_1 \wedge \Omega_1 - \eta_2 \wedge \Omega_2 - \eta_3 \wedge \Omega_3$$
or $s_{11} = -1$ and $\sigma = \sigma_{S(\eta')}$ for some $\eta'$ orthogonal to $\eta = \eta_1$, i.e.,
$$\sigma = \sigma_{\eta,t} = -\eta \wedge \Omega + \cos t \, \mathrm{Re}\, \Psi + \sin t \, \mathrm{Im}\, \Psi.$$
Thus in this case we have nearly parallel $\G_2$-structures of different types sharing the same metric: $\sigma_{S(\eta_1)}$ is of type 1, while $\sigma_{\eta,t}$ are of type 2.

Now we can prove the main result of this section.
\begin{Theorem}\label{thmD1}
Let $(M,\sigma)$ be a compact nearly parallel $\G_2$-manifold which is normalized so that $\tau_0 = 4$ and is not a space of constant curvature. Then:
\begin{enumerate}
\item \label{equiv1} $(M,\sigma)$ is of type 1 if and only if $\dim D_1 = 0$.
\item \label{equiv2} $(M,\sigma)$ is of type 2 if and only if $\dim D_1 = 1$.
\item \label{equiv3} $(M,\sigma)$ is of type 3 if and only if $\dim D_1 = 2$.
\end{enumerate}
\end{Theorem}
\proof Suppose that $(M,\sigma)$ is of type 2. Then the holonomy of the cone is $\SU(4)$ or $\Sp(2)\mathbb{Z}_4 \subset \SU(4)$ and the consideration above show that $\sigma$ is $\sigma_t$ from \pref{sigmat} for some $t$. Again by Proposition~\ref{propD1} we have $\dim D_1 \le 1$. By definition, the contact form $\eta$ of the Sasakian structure satisfies $\nabla \eta = \Omega$. On the other hand, $\eta \lrcorner \sigma_t = -\Omega$ and therefore $\nabla \eta = -\frac{1}{4} \tau_0 \eta \lrcorner \sigma_t$. Thus $\eta \in D_1$, $D_1 = \mathrm{span} \{ \eta \}$ and $\dim D_1 = 1$.

Let $(M,\sigma)$ be of type 3. Then $M$ is $3$-Sasakian and the holonomy of the cone is $\Sp(2)$, so $\sigma = \sigma_{S(\eta)}$ for some $\eta \in S^2$. We shall show that $D_1$ is the orthogonal complement of $\eta$ in $\mathrm{span} \{ \eta_1, \eta_2, \eta_3 \}$. Without loss of generality we can assume that $\eta = \eta_1$ (otherwise we shall change the orthonormal frame $\eta_1, \eta_2, \eta_3$). Then $\sigma \in \{ \sigma_{\eta_2,t} \}$ and $\sigma \in \{ \sigma_{\eta_3,t} \}$, so as above $\eta_2, \eta_3 \in D_1$ and therefore $\dim D_1 \ge 2$. By Proposition~\ref{propD1} every element of $D_1$ induces a Sasakian structure on $(M,g)$ and by Lemma~8.1.17 in \cite{BG} it lies in $\mathrm{span} \{ \eta_1, \eta_2, \eta_3 \}$. Thus, if we assume that $\dim D_1 \ge 3$, we must have $D_1 = \mathrm{span} \{ \eta_1, \eta_2, \eta_3 \}$. But $\nabla \eta_1 = \Omega_1$, while $S(\eta_1)$ is the diagonal matrix with diagonal elements $1,-1,-1$ and $$-\frac{1}{4} \tau_0 \eta_1 \lrcorner \sigma_{S(\eta_1)} = -\eta_1 \lrcorner (\eta_1 \wedge \Omega_1 - \eta_2 \wedge \Omega_2 - \eta_3 \wedge \Omega_3) = -\Omega_1 - \eta_2 \wedge \eta_1 \lrcorner\Omega_2 - \eta_3 \wedge \eta_1 \lrcorner\Omega_3$$
$$= -\Omega_1 - 2\eta_2 \wedge \eta_3 \neq \Omega_1.$$
Hence $\eta_1 \not \in D_1$ and we have $D_1 = \mathrm{span} \{ \eta_2, \eta_3 \} = \eta_1^\perp$ and $\dim D_1 = 2$. This proves \pref{equiv3} since the reverse implication follows from Proposition~\ref{propD1}.

Suppose now that $\dim D_1 = 1$. Then, by Proposition~\ref{propD1}, $(M,g)$ is Sasaki-Einstein but not $3$-Sasakian. Thus, to prove the reverse implication of \pref{equiv2} we only have to show that the case $Hol(\hat M) = \Sp(2)\mathbb{Z}_4$ with $\sigma = \sigma_{S(\eta_1)}$ is impossible. Indeed, the proof of Proposition~\ref{propD1} yields that the contact form $\eta_1$ must be an infinitesimal deformation of $\sigma_{S(\eta_1)}$ but in the type 3 case above we saw that this is not true. This completes the proof of \pref{equiv2}.

Finally, \pref{equiv1} follows from Proposition~\ref{propD1}, \pref{equiv2} and \pref{equiv3}. \qed

\medskip

\begin{Remark}
{\rm A nearly parallel $\G_2$-structure of type 2 is a part of a whole curve $\sigma_t$ of such structures. It is easy to see that $\frac{d\sigma_t}{dt} = *(\eta \wedge \sigma_t)$ ($\eta$ is the contact form of the Sasaki-Einstein structure). Let $\xi$ be the vector field dual to $\eta$. Since all $\sigma_t$ have $\tau_0 =4$ and $D_1 = \mathrm{span} \{ \eta \}$, we obtain from Proposition~\ref{c2} and Proposition~\ref{np} that
$$L_\xi \sigma_t = d(\xi \lrcorner \sigma_t) + \xi \lrcorner d\sigma_t = -\frac{1}{2}d^2 \eta + 4\xi \lrcorner * \sigma_t = -4*(\eta \wedge \sigma_t).$$
Let $\varphi_s$ be the flow of $\xi$. Since $\xi$ is Killing, $\varphi_s$ preserves the metric and thus also $\eta$ and $*$. Now the fact $L_\xi \sigma_t = - \left. \frac{d\varphi_s(\sigma_t)}{ds} \right|_{s=0}$ and the above equations imply
$$\frac{d\varphi_s(\sigma_t)}{ds} = 4*(\eta \wedge \varphi_s(\sigma_t)) = \frac{d\sigma_{t+4s}}{ds}.$$
This and $\varphi_0(\sigma_t) = \sigma_t$ show that $\varphi_s(\sigma_t) = \sigma_{t+4s}$ for all $s$. Thus the flow of $\xi$ acts transitively on the family $\{ \sigma_t \}$ and so the members of this family are equivalent $\G_2$-structures.

In a similar way, if the type is 3, one can generate the whole $D_1$ through curves in the $\mathbb{R}P^2$-family $\{ \sigma_{S(\eta)} \}$. But this family consists of equivalent $\G_2$-structures since a $3$-Sasakian manifold admits an isometric $\SO(3)$ or $\Sp(1)$ action which is transitive on the oriented orthonormal frames $(\eta_1,\eta_2,\eta_3)$ and therefore transitive also on the family $\{ \sigma_{S(\eta)} \}$.

Thus, whatever the type of the nearly parallel $\G_2$-structure, the "interesting" infinitesimal deformations are in the space $D_3$. }
\end{Remark}

\begin{Remark}
{\rm We have seen above that if the holonomy group of the cone $\hat M$ is $\Sp(2)\mathbb{Z}_4$, then $M$ has nearly parallel $\G_2$-structures of different type sharing the same metric and orientation. This is possible because they induce different spin structures on $M$ and therefore also on $\hat M$. Indeed, $\Sp(2)\mathbb{Z}_4$ has two different embeddings in $\Spin(8)$. The first one, $i_1$, is the restriction on $\Sp(2)\mathbb{Z}_4$ of the embedding of $\SU(4)$ in $\Spin(8)$. The second, $i_2$, is equal to $i_1$ on the identity component of $\Sp(2)\mathbb{Z}_4$ and to $-i_1$ on the other component, i.e.,
$$i_2([a,1]) = i_1([a,1]), \quad i_2([a,i]) = -i_1([a,i]) \quad \mbox{for } a \in \Sp(2).$$
Let $E \cong \mathbb{C}^4$ be the standard representation of $\Sp(2)$. Then the spin representation, restricted to $\Sp(2)$, is isomorphic to $\sum_{p=0}^4 \Lambda^p E$. The action of $i_1 (\Sp(2)\mathbb{Z}_4)$ is given by
$$i_1([a,z])\alpha = z^p a\alpha \quad \mbox{for } \alpha \in \Lambda^p E$$
and the space of invariant spinors is $2$-dimensional: $\Lambda^0 E \oplus \Lambda^4 E$. On the other hand, the action of $i_2 (\Sp(2)\mathbb{Z}_4)$ is
$$i_2([a,1])\alpha = a\alpha, \quad i_2([a,i])\alpha = -i^p a\alpha \quad \mbox{for } \alpha \in \Lambda^p E$$
and the space of invariant spinors is $1$-dimensional: $\mathbb{C}\sigma_E \subset \Lambda^2 E = \mathbb{C}\sigma_E \oplus \Lambda^2_0 E$, where $\sigma_E$ is the $\Sp(2)$-invariant symplectic form. Thus an $8$-dimensional manifold with holonomy group $\Sp(2)\mathbb{Z}_4$ is equipped with two canonical spin structures, one of which carries $N=2$ and the other $N=1$ parallel spinors. Similarly, a $7$-dimensional manifold whose cone has holonomy group $\Sp(2)\mathbb{Z}_4$ has two spin structures, with $N=2$ and $N=1$ real Killing spinors respectively. This adds to the results in \cite{Wang}, where in part 2b of Theorem~4.1 $N=1$ is given as the only possibility, while the group $\Sp(2)\mathbb{Z}_4$ is completely missing in part 3 of Corollary~5.2. Notice that the existence of $7$-dimensional manifolds with cones having holonomy group $\Sp(2)\mathbb{Z}_4$ has been proved in \cite{MS}. }
\end{Remark}
%
\section{The $\G_2$-Laplace operator}\label{laplace}
In Section~\ref{sec3} we have seen that infinitesimal deformations of nearly
parallel $\G_2$-manifolds give rise to coclosed eigenforms of the Hodge-de Rham Laplacian acting on sections of $\Lambda^3_{27}\T^* M$. By the classical
Weitzenb\"ock formula the Hodge-de Rham Laplacian is written as
\begin{equation}\label{delta}
\Delta = d^*d + d d^* = \nabla^*\nabla + q(R),
\end{equation}
where $q(R)$ is an endomorphism of the form bundle, which is linear in the curvature $R$ and satisfies
$q(R) = \Ric$ on the space of 1-forms.
We will define the operator $q(R)$ in the following more general setting. Let $(M,g)$ be an $n$-dimensional Riemannian manifold. For a representation $V$ of $\O(n)$ let $VM$ denote the corresponding associated vector bundle. We denote the action of $\alpha \in \Lambda^2 \T^* \cong \mathfrak{so}(n)$ on $V$ by $\alpha_*$ (here $\T$ denotes the standard representation $\mathbb{R}^n$ of $\O(n)$) and in a similar way the action of $\alpha \in \Lambda^2 \T^*_x M$ on $V_x M$, $x \in M$. The endomorphism $q(A) \in \End(VM)$ is defined 
for any $A \in \Lambda^2 \T^* M \otimes \End(VM)$ by
\begin{equation}\label{qrdef}
q(A) = \sum_{i<j}\,(e_i \wedge e_j)_* \,A(e_i \wedge e_j)  ,
\end{equation}
where $\{e_i\}$ is a local orthonormal frame of $\TM$.
Notice that in this definition $\{e_i\wedge e_j \,\, | \,\, i <j \, \}$ could be replaced by any other orthonormal basis of $\Lambda^2 \TM$. The curvature $R$ of the Levi-Civita connection $\nabla$ or, more generally, the curvature $\bar R$ of any metric connection $\bar \nabla$ on $(M,g)$ defines a section of $\Lambda^2 \T^* M \otimes \End(VM)$, thus
the endomorphisms $q(R)$ and $q(\bar R)$ are well defined. We denote by $\bar\Delta $ the Laplace type operator
\begin{equation}\label{bardelta}
\bar \Delta := \bar\nabla^*\bar\nabla + q(\bar R)  .
\end{equation}
The operator $\Delta_L := \nabla^* \nabla + q(R)$ for the Levi-Civita connection $\nabla$ and a sub\-representation $V\subset \otimes^p\,\T$ is also called {\it Lichnerowicz Laplacian} (cf. \cite{besse}, Chapter~1 I). Because of~(\ref{delta}) it coincides
on differential forms with the Hodge deRham Laplacian $\Delta$.

Now let us return to the case of nearly parallel $\G_2$-manifolds. We will call the operator $\bar\Delta$, defined with the canonical connection $\bar\nabla$, 
{\it the $\G_2$-Laplace operator}. In order to compute the spectrum of the Lichnerowicz Laplacian $\Delta_L$ on naturally reductive spaces, it turns out
to be convenient to express $\Delta_L$ through $\bar\Delta$. Thus, our next aim will be to compute the difference $\bar \Delta - \Delta_L$, which we do by calculating the differences
$\bar\nabla^*\bar\nabla - \nabla^*\nabla$ and $q(\bar R) - q(R)$
separately.
A direct calculation using \pref{nablabar2} and the third equation of Lemma~\ref{P2} gives
$$
\begin{array}{rl}
\bar R_{X,Y}Z - R_{X,Y}Z & = (\tfrac{\tau_0}{12})^2 \,[ 2P(P(X,Y),Z) + P(P(Y,Z),X) + P(P(Z,X),Y) ] \\[1.5ex]
&= (\tfrac{\tau_0}{12})^2 \, [ 4P(P(X,Y),Z) - 3g(X,Z)Y + 3g(Y,Z)X ].
\end{array}
$$
Thus\,
$
\bar R (X \wedge Y) - R(X \wedge Y) = (\tfrac{\tau_0}{12})^2 \,[ 4P_{P(X \wedge Y)} - 3(X \wedge Y)_* ].
$
Substituting this equation into the definition of the curvature endomorphisms $q(R)$  and $q(\bar R)$, we obtain
\begin{equation}\label{qr}
q(\bar R) \,-\, q(R)  \;=\; -3 (\tfrac{\tau_0}{12})^2 \,\Cas^{\mathfrak{so}(7)} \,+\;
4(\tfrac{\tau_0}{12})^2 \,S
\end{equation}
where $\Cas^{\mathfrak{so}(n)}$ is the $\so(n)$-Casimir operator  $\sum_{i<j} (e_i \wedge e_j)_*(e_i \wedge e_j)_*$ and $S$ is defined as 
\begin{equation}\label{S}
S = \sum_{i<j} (e_i \wedge e_j)_* P_{P(e_i \wedge e_j)}.
\end{equation}
Since $P: \Lambda^2 \T \cong \Lambda^2_7 \oplus \Lambda^2_{14} \to \T$ is a $\G_2$-equivariant map,  $\left.P\right|_{\Lambda^2_{14}}= 0$ and we may replace in the
sum in (\ref{S}) the orthonormal basis $\{ e_i\wedge e_j \,\, | \,\, i<j \}$ of $\Lambda^2 \T$ with the orthonormal basis $\{f_i=\tfrac{1}{\sqrt 3}e_i \lrcorner \,\sigma \,\, | \,\, i=1,\dots,7\}$ of $\Lambda^2_7$. Because obviously ${f_i}_* = \tfrac{1}{\sqrt 3}(e_i \lrcorner \,\sigma)_* = \tfrac{1}{\sqrt 3}P_{e_i}$ and, by \pref{ai6},
$P(f_i) = \sqrt{3}\,  e_i$, we obtain
$
S = \sum {f_i}_* P_{P(f_i)} = \sum P_{e_i} P_{e_i}  .
$
For the difference $\bar \nabla^*\bar\nabla -\nabla^*\nabla$ of the two rough Laplacians we derive
directly from \pref{nablabar2}
$$
\bar \nabla^*\bar\nabla -\nabla^*\nabla = \sum \left(
\tfrac{\tau_0}{6}P_{e_i}\bar \nabla_{e_i} +
(\tfrac{\tau_0}{12})^2P_{e_i}P_{e_i} \right).
$$
Summarizing these calculation we obtain an expression for the difference of $\bar\Delta$
and $\Delta_L$.
\begin{Proposition}\label{pdiff}
The difference of the Laplace type operators $\bar \Delta$ and $\Delta_L$ on a nearly parallel $\G_2$-manifold is given by
\begin{equation}\label{diff}
\bar \Delta \;-\; \Delta_L
\;\;=\;\;
\tfrac{\tau_0}{6} \sum P_{e_i}\bar \nabla_{e_i} \;-\;\; 3(\tfrac{\tau_0}{12})^2
\Cas^{\mathfrak{so}(7)} \;\;+\;\; 5(\tfrac{\tau_0}{12})^2 \sum P_{e_i}P_{e_i}  .
\end{equation}
\end{Proposition}

We shall apply this result for the space $\Lambda^3_{27}$. Recalling that the $\so(n)$-Casimir operator  acts as $-p(n-p)\id$ on the space of p-forms, we obtain $\Cas^{\mathfrak{so}(7)} \gamma = -12 \gamma$ for $\gamma \in \Lambda^3_{27}$. A straightforward computation on an explicit element (e.g. the element from \pref{expl}) shows that the $\G_2$-equivariant map $\sum P_{e_i}P_{e_i}$ acts as $-8\id$ on $\Lambda^3_{27}$.
Thus it remains to compute $\sum P_{e_i}\bar \nabla_{e_i}$.
The map $\sum P_{e_i} \circ e_i \lrcorner : \T^* \otimes \Lambda^3_{27} \to \Lambda^3$ is $\G_2$-equivariant. Hence, because of \pref{split} and \pref{dec2}, it can be non-zero only on the components of $\T^* \otimes \Lambda^3_{27}$ which are isomorphic to $V_{1,0} \cong \T$ and $V_{2,0} \cong \Lambda^3_{27}$. A straightforward computation on explicit elements shows that
$\sum P_{e_i} \circ e_i \lrcorner$ is $-3* \circ \, \varepsilon$ on the component $\T \subset \T^* \otimes \Lambda^3_{27}$ and $* \circ \varepsilon$ on the component $\Lambda^3_{27} \subset \T^* \otimes \Lambda^3_{27}$. Here again $\varepsilon: \T^* \otimes \Lambda^3 \to \Lambda^4$ denotes the wedging map. Thus
$$\sum P_{e_i}\bar \nabla_{e_i} \gamma = \sum P_{e_i} \circ e_i \lrcorner \bar \nabla \gamma = -3* \varepsilon (\bar \nabla \gamma)_\T + * \varepsilon (\bar \nabla \gamma)_{\Lambda^3_{27}} = -3* (\bar d \gamma)_{\Lambda^4_7} + * (\bar d \gamma)_{\Lambda^4_{27}},$$
where $\bar d := \varepsilon \circ \bar \nabla$. Now, by \pref{nablabar2} we have $\bar d = d - \frac{\tau_0}{12}\sum e_i \wedge P_{e_i}$. Again a simple computation on an element of $\Lambda^3_{27}$ shows that $\sum e_i \wedge P_{e_i} = -2*$ on $\Lambda^3_{27}$ and we
eventually obtain
\begin{Lemma}\label{a3}
On sections of $\Lambda^3_{27}\T^* M$ it holds that
$\;\bar d = d + \tfrac{\tau_0}{6}\ast $. In particular, we have
$$
(\bar d \gamma)_{\Lambda^4_7}
= (d \gamma)_{\Lambda^4_7}, \quad (\bar d \gamma)_{\Lambda^4_{27}} = (d \gamma)_{\Lambda^4_{27}} + \tfrac{\tau_0}{6}\ast \gamma, \quad (d \gamma)_{\Lambda^4_1} = (\bar d \gamma)_{\Lambda^4_1} = 0 \quad \mbox{for }\; \gamma \in \Omega^3_{27}(M).
$$
\end{Lemma}

\noindent
Using this lemma we obtain
$\sum P_{e_i}\bar \nabla_{e_i} \gamma = -3* (d \gamma)_{\Lambda^4_7} + * (d \gamma)_{\Lambda^4_{27}} + \tfrac{\tau_0}{6} \gamma$, which finally enables us to compute
the difference  $\bar \Delta - \Delta$ on $\Omega^3_{27}(M)$.
Combining the formulas above we find

\begin{Proposition}\label{gesamt}
Let $(M^7,g,\sigma)$ be a nearly parallel $\G_2$-manifold and let
$\gamma$ be a 3-form in $\Omega^3_{27}(M)$. Then
$$
\bar \Delta \gamma =  \Delta \gamma  -\tfrac{\tau_0}{2} * (d \gamma)_{\Lambda^4_7} + \tfrac{\tau_0}{6}* (d \gamma)_{\Lambda^4_{27}}  .
$$
In particular, $\Delta$ and $\bar \Delta$ coincide on closed forms in $\Omega^3_{27}(M)$. Moreover, if
$\gamma$ is a 3-form in $\Omega^3_{27}(M)$ with $(d^* \gamma)_{\Lambda^2_7} = 0$, then
$
\bar \Delta\gamma = \Delta \gamma +\tfrac{\tau_0}{6} \ast d\gamma  .
$
\end{Proposition}
\proof
It only remains to prove the last statement:
Recall from the proof of Lemma~\ref{lemf3}, that for $\gamma \in \Omega^3_{27}(M)$ the condition
$(d^* \gamma)_{\Lambda^2_7} = 0$ is equivalent to $(d \gamma)_{\Lambda^4_7} = 0$
and thus, by Lemma~\ref{a3}, also to $(d \gamma)_{\Lambda^4_{27}} = d\gamma$.
Substituting this into the equation  for $\bar \Delta$ implies the last statement.
\qed


\bigskip

%
\section{Infinitesimal Einstein deformations}\label{sec5}
Nearly parallel $\G_2$-structures induce Einstein metrics and thus
infinitesimal deformations of such structures are related to
infinitesimal Einstein deformations. In this section we shall
consider the space of infinitesimal Einstein deformations of a given
nearly parallel $\G_2$-metric and realize it as a direct sum of
certain spaces of 3-forms in $\Omega^3_{27}(M)$.

Let $g$ be an Einstein metric with $\Ric = E g$. From \cite{besse},
Theorem 12.30, the space of infinitesimal Einstein deformations of
$g$ is isomorphic to the set of trace-free symmetric bilinear forms
$h$ on $\TM$ with $\delta h = 0$ and $\Delta_L h = 2E h$, where $\Delta_L
= \nabla^*\nabla + q(R)$ is the so-called Lichnerowicz Laplacian
(see the previous section). Note that for a nearly parallel $\G_2$-metric
the eigenvalue can be written as $2E =\frac{2\scal}{7}=\frac{3\tau^2_0}{4}$.

As a $\G_2$-representation the space $S^2_0\T^*$ is isomorphic to
$\Lambda^3_{27}\T$. We shall now use the explicit identification ${\bf
i}$ in order to identify infinitesimal Einstein deformations with
certain eigenforms of the Laplacian on forms in $\Lambda^3_{27}\T$.
To do this we still need an analogue of Proposition~\ref{gesamt}.

We apply the results of Proposition~\ref{pdiff} to
the space $S^2_0 \T^*$. It is well known that the $\so(n)$-Casimir operator acts on
$S^2_0 \T^*$ as $-2n \Id$, i.e., as $-14 \Id$ in our case. Moreover it is clear that similarly
$\sum P_{e_i} P_{e_i} $, as a $\G_2$-equivariant map,  acts as a
multiple of the identity. An explicit calculation, e.g. on the
element $e^1 \odot e^2$, shows that $\sum P_{e_i} P_{e_i}  = -14 \Id$. 

It remains to determine $\sum P_{e_i}\bar \nabla_{e_i} h$, i.e., $Q(\bar \nabla h)$,
where $Q : \T^* \otimes S^2_0\T^* \rightarrow S^2_0\T^*  $ is the
$\G_2$-equivariant map defined as $Q = \sum P_{e_i} \circ e_i \lrcorner$. The map $Q$ is different from zero only on the component of $T^* \otimes S^2_0 \T^*$, which is isomorphic to $S^2_0 \T^*$. Let
$i_2 : S^2_0 \T^* \rightarrow \T^* \otimes S^2_0 \T^*$ be the embedding
given as $i_2(h) = (1 \otimes \pi_0) \circ C (g \otimes h)$, where
$g$ is the metric, $C : {\T^*}^{\otimes 4} \rightarrow {\T^*}^{\otimes 3}$ is
defined by $C(a\otimes b \otimes c \otimes d )= a \otimes P(b,c)
\otimes d$ and $\pi_0: \T^* \otimes \T^* \rightarrow S^2_0\T^*$ denotes the
standard projection. Moreover, let $\pi_2 : \T^* \otimes S^2_0\T^*
\rightarrow S^2_0\T^*$ be the projection "inverse" to $i_2$, i.e. $\pi_2 \circ
i_2 = \id$ and $\pi_2$ vanishes on the components of $\T^* \otimes S^2_0\T^*$ that are not isomorphic to $S^2_0\T^*$. Then an explicit calculation, e.g. on $e_1 \odot e_2$, shows
that $Q \circ i_2 = -7 \id$ and thus $Q = -7 \pi_2$. Substituting
this and the results for $\Cas^{\mathfrak{so}(7)}$ and $\sum P_{e_i} P_{e_i}$ into equation (\ref{diff}), we obtain
\begin{equation}\label{sym-lap}
(\bar \Delta - \Delta_L) h \;=\; -\tfrac{7\tau_0}{6}\,\pi_2(\bar
\nabla h) \,-\,\tfrac{7\tau^2_0}{36}\, h  .
\end{equation}

Since $S^2_0\T^*$ and $\Lambda^3_{27}$ are isomorphic representations of $\G_2$ and $\bar \nabla$ is a $\G_2$-connection, the bundles $S^2_0\T^*M$ and $\Lambda^3_{27}M$ share the same $\G_2$-Laplace operator $\bar \Delta$, i.e., with the $\G_2$-equivariant isomorphism ${\bf i}: S^2_0\T^* \to \Lambda^3_{27}$ we have ${\bf i} \circ \bar \Delta \circ {\bf i}^{-1} = \bar \Delta$. Hence, to compute ${\bf i} \circ \Delta_L \circ {\bf i}^{-1}$
we need to compute ${\bf i} \circ \pi_2 \circ \bar \nabla \circ {\bf i}^{-1}$. An easy calculation shows
that ${\bf i} = \pi_1 \circ (1\otimes {\bf i}) \circ i_2$, where
$\pi_1 : \T^* \otimes \Lambda^3_{27} \rightarrow \Lambda^3_{27} $ is defined as $\pi_1(\alpha \otimes \gamma)= \frac27\ast(\alpha
\wedge \gamma)_{\Lambda^4_{27}}$. The map $i_2\circ \pi_2$ is the
projection on the component isomorphic to $S^2_0\T^*$ in $\T^* \otimes S^2_0\T^*$
and since $\pi_1 \circ (1\otimes {\bf i})$ is invariant with values
in $S^2_0\T^*$, it vanishes on all other components of $\T^* \otimes
S^2_0\T^*$. Hence $\pi_1 \circ (1\otimes {\bf i}) = {\bf i}\circ
\pi_2$ and we obtain from Lemma~\ref{a3}
\begin{equation}\label{projection}
{\bf i}\circ \pi_2(\bar \nabla h)  = \pi_1 \circ (1\otimes {\bf i})
\bar\nabla h = \pi_1 \bar\nabla ({\bf i} (h)) = \tfrac27\ast (\bar d
{\bf i} (h))_{\Lambda^4_{27}} = \tfrac27\ast (d {\bf i} (h))_{\Lambda^4_{27}} +
\tfrac{\tau_0}{21}{\bf i} (h)  .
\end{equation}

Let $h$ be an infinitesimal Einstein deformation and let $\gamma \in \Omega^3_{27}(M)$ denote the 3-form ${\bf i}(h)$. Then the condition $\delta h =0$ translates into
$(d^* \gamma)_{\Lambda^2_7} = 0$ or, equivalently, to $d\gamma= (d\gamma)_{\Lambda^4_{27}}$. Indeed, by Remark~\ref{r1} we have that $(\nabla h)_\T = (\bar \nabla h)_\T$ and $(\nabla \gamma)_\T = (\bar \nabla \gamma)_\T$. Thus $\delta h = 0$ is equivalent to $(\nabla h)_\T = 0$ and also to $(\bar\nabla h)_\T
= 0$. But since ${\bf i}$ is an $\G_2$-equivariant map, $(\bar\nabla h)_\T = 0$ if and only if $(\bar\nabla \gamma)_\T = 0$, i.e., $(\nabla \gamma)_\T = 0$. However this is equivalent to $(d^* \gamma)_{\Lambda^2_7} = 0$ and also to $(d\gamma)_{\Lambda^4_7} = 0$. Then by Lemma~\ref{a3} $(d\gamma)_{\Lambda^4_7} = 0$ can be written as $d\gamma= (d\gamma)_{\Lambda^4_{27}}$.

Finally, we apply ${\bf i}$ to (\ref{sym-lap}), use Proposition~\ref{gesamt} and substitute
(\ref{projection}) to obtain

\begin{Proposition}
For each $\gamma \in \Omega^3_{27}(M)$ the following equation is satisfied:
\begin{equation}\label{sym-lap2}
{\bf i} \; \Delta_L \, {\bf i}^{-1} (\gamma)
\;=\;
\Delta \gamma  \;-\;  \tfrac{\tau_0}{2} * (d \gamma)_{\Lambda^4_7}  \;+\;
\tfrac{\tau_0}{2}* (d \gamma)_{\Lambda^4_{27}}  \;+\;  \tfrac{\tau_0^2}{4}\gamma  .
\end{equation}
\end{Proposition}

With this formula we are able to translate the conditions for
infinitesimal Einstein deformations into equivalent conditions for
$3$-forms in $\Omega^3_{27}(M)$: The traceless symmetric bilinear form $h$
is an infinitesimal Einstein deformation if and only if
$\gamma = {\bf i}(h)$ is a section of $\Lambda^3_{27}\T^* M$ with
$(d\gamma)_{\Lambda^4_7}=0$ (or, equivalently,
$(d^*\gamma)_{\Lambda^2_7}=0$), satisfying the equation
\begin{equation}\label{inf-ein}
\Delta \gamma \;+\; \tfrac{\tau_0}{2}\ast d\gamma  \;-\; \tfrac{\tau_0^2}{2}\gamma
\;=\; 0 .
\end{equation}
We want to decompose the solution space of this equation into
eigenspaces of the operator $\ast d$. This is possible since $\ast
d$ is a symmetric operator, commuting with the operator on the left
hand side of equation \pref{inf-ein} and preserving the condition
$(d^*\gamma)_{\Lambda^2_7}=0$. Indeed, $\ast d \alpha$ is coclosed for any
differential form $\alpha$. Moreover, the solution space
is finite dimensional because it is the kernel of an elliptic operator.
Assume that
$\ast d \gamma = \lambda \gamma$ with $\lambda \neq 0$. Then $\gamma$ is
coclosed (in particular, $(d^*\gamma)_{\Lambda^2_7}=0$ and $(d\gamma)_{\Lambda^4_7}=0$) and \pref{inf-ein} yields the
quadratic equation
\begin{equation}\label{inf-ein1}
\lambda^2  \;+\; \tfrac{\tau_0}{2}\lambda  \;-\; \tfrac{\tau_0^2}{2} \;=\; 0
\end{equation}
with the solutions $\lambda= -\tau_0$ and
$\lambda=\frac{\tau_0}{2}$. In the case $\lambda = 0$ we obtain
$d\gamma = 0$ and $dd^*\gamma = \frac{\tau^2_0}{2}\gamma$. Moreover a
solution $\gamma$ of the last equation is automatically closed and thus
$(d\gamma)_{\Lambda^4_7}=0$ as well as  $(d^*\gamma)_{\Lambda^2_7}=0$. Summarizing
we have

\begin{Theorem}\label{ges-ein}
Let $(M, \,\sigma,\, g)$ be a compact nearly parallel $\G_2$-manifold. Then the
space of infinitesimal Einstein deformations of $g$ is isomorphic to
the direct sum of the spaces
$$
\{ \gamma \in \Omega^3_{27} \,\, | \,\, *d\gamma = -\tau_0 \gamma \}, \quad
\{ \gamma \in \Omega^3_{27} \,\, | \,\, *d\gamma = \tfrac{\tau_0}{2} \gamma \}, \quad \{ \gamma \in \Omega^3_{27} \,\, | \,\, dd^* \gamma = \tfrac{\tau_0^2}{2} \gamma \}.
$$
\end{Theorem}

Notice that the first space is the space $D_3$ from Theorem~\ref{thmid}. Thus any element $f_3 \in D_3$ satisfies ${\bf i} \; \Delta_L {\bf i}^{-1} \,(f_3) = \frac{3\tau_0^2}{4} \,
f_3$, or, equivalently, $\Delta_L\, {\bf j} (f_3) = \frac{3\tau_0^2}{4} \,{\bf j} (f_3)$, which finishes the proof of Theorem~\ref{thmid}.

In order to check in the examples whether or not infinitesimal
Einstein deformations exist it will be convenient to embed these
three spaces into eigenspaces of the operator $\bar\Delta$
acting on sections of $\Lambda^3_{27}\T^* M$. Let $\gamma$ be a 3-form as
above with $\ast d \gamma = \lambda \gamma$ for $\lambda \neq 0$. Then
$\gamma$ is coclosed and $\Delta \gamma = \lambda^2 \gamma$. Thus
Proposition~\ref{gesamt} implies $\bar \Delta \gamma = (\lambda^2 +
\frac{\tau_0}{6} \lambda) \gamma$. In the case $dd^*\gamma =
\frac{\tau^2_0}{2}\gamma$ it follows that $\gamma$ is closed and
we obtain $\bar \Delta \gamma = \Delta \gamma = dd^*\gamma =
\frac{\tau_0^2}{2}\gamma$. This proves

\begin{Lemma}\label{eigenv}
The three summands of
Theorem~\ref{ges-ein} are contained in the eigenspaces of
$\bar\Delta$ acting on 
$\Omega^3_{27}(M)$ for the eigenvalues
$\frac{5\tau^2_0}{6}$, $\frac{\tau^2_0}{3}$ and $\frac{\tau^2_0}{2}$ respectively.
\end{Lemma}

%
%
\section{Naturally reductive spaces}\label{sec-nat-red}
In this section we will make some general remarks which will help us to compute the infinitesimal Einstein deformations of nearly parallel $\G_2$-manifolds that are
naturally reductive homogenous spaces, i.e. reductive spaces where the torsion
of the canonical homogenous connection can be considered as a 3-form. 

\begin{Lemma}\label{hlem}
Let $G/H$ be a $7$-dimensional oriented naturally reductive homogeneous space with
reductive decomposition $\mathfrak{g} = \mathfrak{h} \oplus \mathfrak{m}$.
Suppose that at the initial point $o$ the torsion of the canonical homogeneous
connection $\hat \nabla$
is $\hat
\T_o = -\frac{\tau_0}{6} \sigma_o$ with $\tau_0 \neq 0$ and that $\sigma_o$ is stable and induces the given metric and orientation on $\mathfrak{m}$. Then $\sigma_o$ defines by translations a $G$-invariant $3$-form $\sigma$ and thus a $\G_2$-structure on $G/H$ compatible with
the given metric and orientation. This $\G_2$-structure is nearly parallel and its canonical connection is $\bar \nabla = \hat \nabla$. In particular, $d \sigma = \tau_0 *\sigma$.

Moreover if $G/H$ is standard up to a scaling factor $c^2$, i.e., $\mathfrak{m}$ is the orthogonal complement of $\mathfrak{h}$ with respect to the Killing form $B$ of
$\mathfrak{g} $ and the metric is induced by the restriction of $-c^2 B$ to $\mathfrak{m}$, then the scalar curvature is $\scal = \frac{63}{20c^2}$ and $\tau_0^2 = \frac{6}{5c^2}$.
\end{Lemma}

\proof
Since $\sigma_o = -\frac{6}{\tau_0} \hat \T_o$ is an $H$-invariant $3$-form on $\mathfrak{m}$, $\sigma = -\frac{6}{\tau_0} \hat \T$ is a $G$-invariant $3$-form on $G/H$. In particular, $\sigma$ is parallel with respect to the canonical homogeneous connection
$$
\hat \nabla = \nabla + \tfrac{1}{2} \,\hat \T = \nabla - \tfrac{\tau_0}{12}\,\sigma  .
$$

For $X \in \mathbb{R}^7$ we have the identity $P_X \sigma = 3X \lrcorner *\sigma$ (which follows from \pref{ai3} and \pref{ai4} or by an explicit computation for some $X \neq 0$).
Thus
$$
\nabla_X \sigma = \hat \nabla_X \sigma + \tfrac{\tau_0}{12}\, P_X \sigma =
\tfrac{\tau_0}{4}X \,\lrcorner *\sigma
$$
and Proposition~\ref{np} implies that the $\G_2$-structure is nearly parallel with $d \sigma = \tau_0 *\sigma$. Moreover $\hat \nabla$ coincides with the canonical connection $\bar \nabla$ of the $\G_2$-structure because of~(\ref{nablabar2}).

Suppose now that $G/H$ is standard (up to a scaling factor $c^2$). Obviously, it is enough to prove the statement about the scalar curvature when $c=1$. Recall that $\hat \T_o(X,Y) = - [X,Y]_\mathfrak{m}$ for $X,Y \in \mathfrak{m}$. Then, considering again $\hat \T_o$ as a $3$-form and using (7.39) in \cite{besse}, we obtain
$$
\scal \;=\; -\tfrac{6}{4} \,|\hat \T_o|^2 \,+\; \tfrac{7}{2}
\;=\; -\tfrac{\tau_0^2}{24}\, |\sigma_o|^2 \,+\; \tfrac{7}{2} .
$$
Since $\sigma_o$ induces the metric on $\mathfrak{m}$ we have $|\sigma_o|^2 = 7$
and, using \pref{scal} to replace $\scal$, the equation above yields $\tau_0^2 = \frac{6}{5}$ and therefore $\scal = \frac{21}{8}\tau_0^2 = \frac{63}{20}$.
\qed

In view of this lemma and the results of the previous section it will be useful to have an algebraic description of some differential operators on naturally reductive spaces.

Let $M = G/H$ be a reductive homogeneous space and $\rho$ be a representation of $H$ on a vector space $V$. Denote by $E:= G \times_\rho V$ the associated vector bundle over $M$. If a $G$-invariant metric is fixed on $M$, then the canonical homogeneous connection $\hat \nabla$ is a metric connection and, as explained in Section~\ref{laplace}, we can define the Laplace type operator $\hat \Delta_\rho = \hat \nabla^* \hat \nabla + q(\hat R)$
acting on sections of $E$. With the same proof as for Lemma~5.2 in \cite{andrei-uwe}
we have

\begin{Lemma}\label{spectrum}
Let $G$ be a compact semi-simple Lie group,
$H \subset G$ a compact subgroup and let $M=G/H$ be standard (up to a scaling factor $c^2$). Then the endomorphism
$q(\hat R)$ acts fibrewise on $E$ as $- \frac{1}{c^2}\,\Cas^H_\rho$ and the operator
$\hat \Delta_\rho$ acts on $\Gamma( E)$, considered as a $G$-representation via the left-regular representation
$l$, as $- \frac{1}{c^2} \,\Cas^G_l$, where the Casimir operator $\Cas^G_{V_\gamma}$ of a $G$-representation
$V_\gamma$ is defined with respect to the Killing form of $G$.
\end{Lemma}

Lemma~\ref{spectrum} can be used to compute the spectrum of $\hat \Delta_\rho$.
We recall that the Peter-Weyl theorem and the Frobenius
reciprocity yield the following decomposition of the left-regular
representation of $G$ into irreducible summands:
\begin{equation}\label{frob}
\Gamma(E) \;\cong\; \overline \bigoplus \, V_\gamma \otimes\Hom_H (V_\gamma, V)  ,
\end{equation}
where the sum is taken over the set of (non-isomorphic) irreducible
$G$-representations $V_\gamma$, labeled by their highest weight $\gamma$.
The Casimir operator acts on $V_\gamma$ as a certain multiple of the identity, which can be
computed explicitly by the Freudenthal formula. Hence the eigenspace of $\hat \Delta_\rho$ for the eigenvalue $\lambda$ is isomorphic as a $G$-representation to the direct sum of the spaces $V_\gamma \otimes \Hom_H (V_\gamma, V)$ for which $\Cas^G_{V_\gamma} = -c^2 \lambda$.


\begin{Corollary}\label{Cas-eig}
Let $G/H$ be standard (up to a scaling factor $c^2$), satisfying the assumptions
of Lemma~\ref{hlem}. Then the eigenspaces of the $\G_2$-Laplace operator
$\bar \Delta$ on $\Omega^3_{27}(M)$ for the eigenvalues $\frac{5\tau^2_0}{6}$, $\frac{\tau^2_0}{3}$, $\frac{\tau^2_0}{2}$ are isomorphic as $G$-representations to the direct sum of spaces $\;V_\gamma \otimes \Hom_H (V_\gamma, \Lambda^3_{27} \mathfrak{m}^*)$, on which the Casimir operator $\Cas^G_{V_\gamma}$ acts as $-1$, $-\frac{2}{5}$, $-\frac{3}{5}$.
\end{Corollary}

\bigskip

In the examples below we have to solve equations of the form $\bar d \varphi + c \ast \varphi = 0$ for 3-forms $\varphi$
on naturally reductive spaces $M=G/H$. Using the explicit embedding of $V_\gamma \otimes\Hom_H (V_\gamma, V)$ of
\pref{frob} into $\Gamma(E)$ we will translate this into an algebraic equation.

As above let $E:= G \times_\rho V$ be the vector bundle over $M=G/H$ associated to
a representation $\rho: H \rightarrow \Aut(V)$.
The space of $H$-equivariant functions from $G$ to $V$, i.e.,
functions $ f:G \to V$ with $f \circ R_h = \rho(h^{-1}) \circ f$ for all  $h \in H$,
can be identified with the space of the sections of $E$. Indeed, the section
$\varphi$  corresponding to the function $f$ is given by $\varphi(\pi(a)) = a(f(a)).$
Here $\pi : G \to G/H$ denotes the projection, $\pi(a) = aH$, and $a \in G$ is considered as a linear isomorphism from $V$ to the fibre $E_{\pi(a)}$, defined on $v\in V$ as
$a(v):= [a,v] \in E_{\pi(a)}$.
Since $G$ acts from the left on the space of $H$-equivariant functions from $G$ to $V$ by
$a \cdot f := L_{a^{-1}}^* f = f \circ L_{a^{-1}},$ we obtain a left action of $G$ on
$\Gamma(E)$.

Let $U$ be an irreducible $G$-representation. Then $U \otimes \Hom_H (U,V)$
embeds into $\Gamma(E)$ by
$$
U \otimes \Hom_H (U,V) \ni \alpha \otimes A
\mapsto f^A_\alpha, \; \mbox{ where } \; f^A_\alpha : G \to V, \quad f^A_\alpha (a) = A(a^{-1} \alpha) .
$$
In particular, fixing $A \in \Hom_H (U,V)$ one obtains a $G$-equivariant
homomorphism $U \to \Gamma(E)$,  given by $U \ni \alpha \mapsto f^A_\alpha$.
The meaning of \pref{frob} is that each $G$-equivariant homomorphism
$U\rightarrow \Gamma(E)$
is obtained in this way. In other words, a subspace of
$\Gamma(E)$ is isomorphic as a $G$-representation to $U$ if and only
if it coincides with the space $\{ f^A_\alpha : \alpha \in U \}$ for
some $A \in \Hom_H (U,V)$, $A \neq 0$.

Let $M = G/H$ be reductive  with $\Ad(H)$-invariant decomposition
$\mathfrak{g} = \mathfrak{h} \oplus \mathfrak{m}$ and $E = \Lambda^s T^*M$, i.e., the vector
bundle associated to the $H$-representation $V = \Lambda^s \m^*$. Then a
straightforward computation shows that in this case $a \cdot \varphi = L_{a^{-1}}^* \varphi$
for $a \in G$ and $\varphi \in \Gamma(E)$. This means that if $\varphi$ corresponds to the
function $f$, then $L_a^* \varphi$ corresponds to the function $L_a^* f$. Let $\hat \nabla$
be the canonical homogeneous connection and consider the operator $\hat d = \varepsilon
\circ \hat \nabla : \Gamma(\Lambda^s T^*M) \to \Gamma(\Lambda^{s+1} T^*M)$. Since $\hat \nabla$
is translation invariant, we have
$$
(\hat d \varphi)_{\pi(a)} = L_{a^{-1}}^* ((\hat d L_a^* \varphi)_{\pi(e)}).
$$
For $(\hat d L_a^* \varphi)_{\pi(e)}$ we obtain the equation
$$
(\hat d L_a^* \varphi)_{\pi(e)}(X_1,\dots,X_{s+1}) \;=\; \sum_{i=1}^{s+1} (-1)^{i-1}
dL_a^* f (X_i)(X_1,\dots,\hat{X}_i,\dots,X_{s+1})
$$
for $X_1,\dots,X_{s+1} \in \mathfrak{m} \cong T_{\pi(e)} M$.
Let $\varphi$ correspond to $f^A_\alpha$. Then
$$dL_a^* f^A_\alpha (X) = (d A(b^{-1}a^{-1} \alpha))_{b=e} (X) = A((d(b^{-1}))_{b=e} (X)
 \cdot a^{-1} \alpha) = -A(X \cdot a^{-1} \alpha),$$
where $X \cdot \alpha$ denotes the action of $X \in \mathfrak{g}$ on $\alpha \in U$. Thus
$$(\hat d L_a^* \varphi)_{\pi(e)}(X_1,\dots,X_{s+1}) = \sum_{i=1}^{s+1} (-1)^i A(X_i
\cdot a^{-1} \alpha)(X_1,\dots,\hat{X}_i,\dots,X_{s+1})$$
for $X_1,\dots,X_{s+1} \in \mathfrak{m} \cong T_{\pi(e)} M$.
In a similar way one obtains
$$
(d L_a^* \varphi)_{\pi(e)}(X_1,\dots,X_{s+1}) = \sum_{i=1}^{s+1} (-1)^i A(X_i \cdot a^{-1}\alpha)(X_1,\dots,\hat{X}_i,\dots,X_{s+1})
$$
$$
+ \sum_{1 \le i < j \le s+1} (-1)^{i+j} A(a^{-1}\alpha)([X_i,X_j]_\mathfrak{m},X_1,\dots,\hat{X}_i,\dots,\hat{X}_j,\dots,X_{s+1})
$$
for $X_1,\dots,X_{s+1} \in \mathfrak{m} \cong T_{\pi(e)} M$.
From these formulas  one can compute $(\hat d \varphi)_{\pi(a)}$ and $(d \varphi)_{\pi(a)}$ for any $a \in G$.

Next we fix a $G$-invariant metric and an orientation on $M$. Then
$\;
*\varphi = L_{a^{-1}}^* ((* L_a^* \varphi)_{\pi(e)}).
$
Therefore, if we would like to solve the $G$-invariant equation
$\;
\hat d \varphi + c *\varphi = 0,
$
for a certain constant $c$, it is enough to solve
$\;
(\hat d L_a^* \varphi)_{\pi(e)} + c* (L_a^* \varphi)_{\pi(e)} = 0 \quad \mbox{for all } a \in G.
$
In fact, we shall be interested in subspaces of solutions of this equation, which are isomorphic
to a given irreducible $G$-representation $U$. Thus we have to find $A \in \Hom_H(U,V)$ so that
$$\sum_{1 \le i_1 < \dots < i_{s+1} \le n} \sum_{j=1}^{s+1} (-1)^i A(e_{i_j} \cdot a^{-1} \alpha)
(e_{i_1},\dots,\hat{e}_{i_j},\dots,e_{i_{s+1}})e^{i_1 \dots i_{s+1}} + c*A(a^{-1}\alpha) = 0$$
for all $a \in G$, $\alpha \in U$. Here $e_1,\dots,e_n$ is a basis of $\mathfrak{m}$. It is clear
that it suffices to write $a=e$ in this equation, i.e., we are looking for $A \in \Hom_H(U,V)$ so that
\begin{equation}\label{main}
\sum_{1 \le i_1 < \dots < i_{s+1} \le n} \sum_{j=1}^{s+1} (-1)^j A(e_{i_j} \cdot \alpha)(e_{i_1},
\dots,\hat{e}_{i_j},\dots,e_{i_{s+1}})\, e^{i_1 \dots i_{s+1}} \;+ \;c*A(\alpha) \;=\; 0
\end{equation}
holds for all $\alpha \in U$. Notice that this equation is $H$-invariant.

%
\section{Examples}
In this section we shall compute the infinitesimal Einstein deformations
of three examples of proper nearly parallel $\G_2$-structures on standard homogeneous
spaces (up to a factor).

The first example is $\SO(5)/\SO(3)$, where the embedding of $\SO(3)$ in $\SO(5)$
is given by the $5$-dimensional irreducible representation of $\SO(3)$. This
space is isotropy irreducible. In fact, the isotropy representation is the
unique $7$-dimensional irreducible representation of $\SO(3)$, which also
defines an embedding of $\SO(3)$ in $\G_2$ and thus a $\G_2$-structure on
$\SO(5)/\SO(3)$. The $\G_2$-structure is proper nearly parallel (cf. \cite{Bryant0}).

The other two examples come from $3$-Sasakian geometry. Recall that there is a second
Einstein metric in the canonical variation of a $3$-Sasakian
metric. In the $7$-dimensional case this metric is induced by a
proper nearly parallel $\G_2$-structure \cite{uwe1}. In general,
for each simply connected compact simple Lie group $G$ there
exists exactly one simply connected $3$-Sasakian homogeneous
manifold of the form $G/H$ and the only other $3$-Sasakian homogeneous
manifolds are the real projective spaces \cite{BG}. The second
Einstein metric is also $G$-homogeneous but not normal
(neither is the $3$-Sasakian metric). But if one writes
the space in the form $\frac{G \times \Sp(1)}{H \times \Sp(1)}$,
it becomes normal \cite{A} and in the $7$-dimensional case even standard (up to a factor).

The simply connected $7$-dimensional homogeneous $3$-Sasakian manifolds are the round
sphere $S^7$ and the Aloff-Wallach space $N(1,1)$. The corresponding second Einstein
metrics are the standard homogenous metrics (up to a factor) on $\frac{\Sp(2) \times
\Sp(1)}{\Sp(1) \times \Sp(1)}$ (the so-called squashed sphere) and  on $\frac{\SU(3)
\times \Sp(1)}{U(1) \times \Sp(1)}$. As remarked by B.~Wilking in \cite{Wilking}, the
latter space was overlooked in the Berger classification of normal homogenous spaces
of positive sectional
curvature. It follows from Equation~(7.87b) of \cite{besse} that a normal homogenous
space has non-negative sectional curvature. However, if in addition one has a
$\G_2$-structure as desribed in Lemma~\ref{hlem} the torsion is non-degenerate and the
sectional curvature has to be positive. Thus by the Berger classification there are
only the examples considered above.

To compute the space of infinitesimal Einstein deformations on our examples $M=G/H$,
we shall proceed in the following way. First we determine which $H$-representation
$V$ define the bundle $\Lambda^3_{27} \T^* M$ and then we use Corollary~\ref{Cas-eig} to
find the irreducible $G$-representations $U$ appearing in the eigenspaces  of $\bar
\Delta$ for the eigenvalues $\frac{5\tau^2_0}{6}$, $\frac{\tau^2_0}{3}$ and $\frac{\tau^2_0}{2}$
(as given in Lemma~\ref{eigenv}). In all three examples the computation of the Casimir
eigenvalues will show that the eigenvalues $\frac{\tau^2_0}{3}$ and $\frac{\tau^2_0}{2}$ do not
appear and thus the spaces of infinitesimal Einstein deformations and infinitesimal $\G_2$-deformations
coincide.  It is interesting to note that in all three cases the non-zero candidates
$U$ comming from the eigenvalue $\frac{5\tau^2_0}{6}$ turn out to be exactly the components
of the adjoint representation of $G$. However this is not too surprising since the Casimir
eigenvalue of the adjoint representation (with respect to the Killing form) is always~$-1$.
If such representations $U$ do exist, we have to solve the equation
$d\varphi = - \tau_0 \ast \varphi$ or equivalently
$\;
(\bar d +  c  \ast) \varphi = 0
\;$
for the constant $c= \frac56\,\tau_0$. By the results of the previous section, this
is reduced to finding $A \in \Hom_H (U,\Lambda^3_{27}\m)$ that solves \pref{main}
with $c= \frac56\,\tau_0$.

For reference below we mention the following facts about Casimir operators. The Casimir operator
of the representation $V(k_1,\dots,k_n)$ of $\Sp(n)$ with highest weight
$\gamma = (k_1,\dots,k_n)$, where $k_1 \ge \dots \ge k_n \ge 0$ are integers, is given by
\begin{equation}\label{Cas-Sp}
\Cas^{\Sp(n)}_{V(k_1,\dots,k_n)} = -\tfrac{1}{4(n+1)} \sum_{i=1}^n \,(2(n-i+1)k_i + k_i^2)
\end{equation}
and the Casimir operator of the representation $V(k_1,\dots,k_n)$ of $\SU(n)$ with highest
weight $\gamma = (k_1,\dots,k_n)$, where $k_1 \ge \dots \ge k_n$ are integers satisfying
$-\frac{n}{2} < k_1 + \dots + k_n \le \frac{n}{2}$, by
\begin{equation}\label{Cas-SU}
\Cas^{\SU(n)}_{V(k_1,\dots,k_n)} = -\tfrac{1}{2n} \, \sum_{i=1}^n ((n+1-2i)k_i + k_i^2) +
\tfrac{1}{2n^2} \, ( \sum_{i=1}^n k_i )^2  .
\end{equation}
Finally, if $V_1$ and $V_2$ are representations of the groups $G_1$ and $\G_2$ respectively, then
\begin{equation}\label{Cas-product}
\Cas^{G_1 \times \G_2}_{V_1 \otimes V_2} = \Cas^{G_1}_{V_1} + \Cas^{ \G_2}_{ V_2}  .
\end{equation}

\bigskip

%
%
\subsection{\large{The example $\SO(5)/\SO(3)$ }} ~ \\

We have the reductive, i.e. $\Ad(\SO(3))$-invariant, decomposition $\mathfrak{so}(5) = \mathfrak{so}(3)
\oplus \mathfrak{m}$, where $\mathfrak{m}$ is the orthogonal complement of $\mathfrak{so}(3)$
with respect to the Killing form of $\mathfrak{so}(5)$. As mentioned above,
$\mathfrak{m}$ is the irreducible $7$-dimensional representation of $\SO(3)$. The complex
irreducible $\SO(3)$-representations can be written as the symmetric powers $S^{2k}E$, where
$E=\C^2$ is the standard representation of the double cover $\Sp(1)$ of $\SO(3)$, in particular
$\mathfrak{m}^\mathbb{C} \cong S^6 E$.
It is easy to obtain the following
decomposition into irreducible summands
$$
{\Lambda^3 \mathfrak{m}^*}^\mathbb{C} \cong \Lambda^3 S^6 E = \C \oplus S^4E \oplus S^6E \oplus S^8E
\oplus S^{12}E  .
$$
We see that there is a 1-dimensional space of $\SO(3)$-invariant
3-forms, which implies that on $M=\SO(5)/\SO(3)$ the canonical
homogeneous connection coincides with the canonical
$\G_2$-connection. Moreover, since $\Lambda^3\m^* \cong \R \oplus \m
\oplus \Lambda^3_{27}\m^*$ as a $\G_2$-representation, we obtain
\begin{equation}\label{deco}
{\Lambda^3_{27}\m^*}^\C \cong S^4E \oplus S^8E \oplus S^{12}E 
.
\end{equation}
Since $\Sp(2)$ double covers $\SO(5)$, the two groups have the same Casimir operator. Therefore,
by Corollary~\ref{Cas-eig} and \pref{Cas-Sp}, we have to find all pairs of integers $(k_1,k_2)$
with $k_1 \ge k_2 \ge 0$, such that $-\frac{1}{12} (4k_1 + k_1^2 + 2k_2 + k_2^2)$ is equal to
one of $-1$, $-\frac{2}{5}$, $-\frac{3}{5}$. The only solution is $(k_1,k_2) = (2,0)$ for the
eigenvalue $-1$. The representation $V(2,0)$ is the adjoint representation of $\Sp(2)$ and it
corresponds, of course, to the adjoint representation $\mathfrak{so}(5)$ of $\SO(5)$.
It remains to
compute the dimension of $\Hom_{\SO(3)}(\mathfrak{so}(5)^\mathbb{C},
{\Lambda^3_{27}\m^*}^\C)$, which turns out to be zero. Indeed, from the reductive decomposition
above we have the following decomposition of $\mathfrak{so}(5)^\mathbb{C}$ into irreducible
$\SO(3)$-representations:
$$
\so(5)^\C \cong \so(3)^\C \oplus \m^\C \cong S^2 E \oplus
S^6 E  .
$$
Comparing this with (\ref{deco}), we see that $\so(5)^\C$ and ${\Lambda^3_{27}\m^*}^\C$
do not have any common components and therefore $\Hom_{\SO(3)}(\mathfrak{so}(5)^\mathbb{C},
{\Lambda^3_{27}\m^*}^\C) = 0$. Thus the eigenvalues $\frac{5\tau^2_0}{6}$, $\frac{\tau^2_0}{3}$
and $\frac{\tau^2_0}{2}$ do not appear in the $\bar \Delta$-spectrum on $\Omega^3_{27}(M)$
and so we have proved

\begin{Proposition}
There are no infinitesimal Einstein deformations and, in particular, no infinitesimal
$\G_2$-deformations of the nearly parallel $\G_2$-structure on $\SO(5)/\SO(3)$.
\end{Proposition}

\bigskip

%
%
\subsection{\large{The example $\frac{\Sp(2) \times \Sp(1)}{\Sp(1) \times \Sp(1)}$ }} ~ \\

\noindent
We denote by $\Sp(1)_u$ and $\Sp(1)_d$ the following embeddings of $\Sp(1)$ in $\Sp(2)
\times \Sp(1)$:
$$
\Sp(1)_u := \{ ( \begin{pmatrix} a & 0 \\ 0 & 1 \end{pmatrix}, 1): a \in \Sp(1) \}, \qquad
\Sp(1)_d := \{ ( \begin{pmatrix} 1 & 0 \\ 0 & a \end{pmatrix}, a): a \in \Sp(1) \}.
$$
In this realization the Lie algebras of $\Sp(1)_u$ and $\Sp(1)_d$ are given as
$$
\mathfrak{sp}(1)_u := \{ ( \begin{pmatrix} a & 0 \\ 0 & 0 \end{pmatrix}, 0): a \in
\mathfrak{sp}(1) \}, \qquad \mathfrak{sp}(1)_d := \{ ( \begin{pmatrix} 0 & 0 \\ 0 & a
\end{pmatrix}, a): a\in \mathfrak{sp}(1) \}  .
$$
We consider the homogeneous space $\frac{\Sp(2) \times \Sp(1)}{\Sp(1)_u \times
\Sp(1)_d}$ as a normal homogeneous space taking the metric induced by $-\frac{1}{24}B$,
where $B$ is the Killing form of $\mathfrak{g} = \mathfrak{sp}(2) \oplus \mathfrak{sp}(1)$.
Then we have the reductive decomposition
$\mathfrak{g} = \mathfrak{h} \oplus \mathfrak{m}$, with
$$
\mathfrak{h} = \mathfrak{sp}(1)_u \oplus \mathfrak{sp}(1)_d, \qquad \mathfrak{m}
= \mathfrak{h}^\bot = \mathfrak{sp}(1)_o \oplus \mathfrak{m}'  .
$$
The Lie algebra $\mathfrak{sp}(1)_o$ and the space $\mathfrak{m}'$ are given as
$$
\mathfrak{sp}(1)_o := \{ ( \begin{pmatrix} 0 & 0 \\ 0 & 2a \end{pmatrix}, -3a): a \in
\mathfrak{sp}(1) \}, \qquad \mathfrak{m}' := \{ ( \begin{pmatrix} 0 & x \\ -\bar{x} & 0
\end{pmatrix}, 0): x \in \mathbb{H} \}  .
$$
We define the orientation by means of the following orthonormal frame of $\mathfrak{m}$:
$$e_1 := \frac{1}{\sqrt{5}}(\begin{pmatrix} 0 & 0 \\ 0 & 2i \end{pmatrix}, -3i), \quad
e_2 := \frac{1}{\sqrt{5}}(\begin{pmatrix} 0 & 0 \\ 0 & 2j \end{pmatrix}, -3j), \quad
e_3 := \frac{1}{\sqrt{5}}(\begin{pmatrix} 0 & 0 \\ 0 & 2k \end{pmatrix}, -3k),$$
$$e_4 := (\begin{pmatrix} 0 & 1 \\ -1 & 0 \end{pmatrix}, 0), \quad
e_5 := (\begin{pmatrix} 0 & i \\ i & 0 \end{pmatrix}, 0), \quad
e_6 := (\begin{pmatrix} 0 & j \\ j & 0 \end{pmatrix}, 0), \quad
e_7 := (\begin{pmatrix} 0 & k \\ k & 0 \end{pmatrix}, 0).$$
Then, computing the comutators of these basis elements, we see that at the initial point $o$
the torsion of the canonical homogeneous connection
is $\hat \T_o = \frac{2}{\sqrt{5}} \sigma_o$, where $\sigma_o$ is given by \pref{sigma}.
Hence, by Lemma~\ref{hlem} we obtain a nearly parallel $\G_2$-structure on $\frac{\Sp(2)
\times \Sp(1)}{\Sp(1)_u \times \Sp(1)_d}$ with $\tau_0 = -\frac{12}{\sqrt{5}}$.

We want to determine infinitesimal Einstein deformations of this structure. Thus
by Corollary~\ref{Cas-eig} together with \pref{Cas-product} and \pref{Cas-Sp}
we are looking for $k_1 \ge k_2 \ge 0$ and $l \ge 0$ such that
$$
\Cas^{\Sp(2) \times \Sp(1)}_{V(k_1,k_2) \otimes V(l)} = -\tfrac{1}{12}\,
(4k_1 + k_1^2 + 2k_2 + k_2^2) - \tfrac{1}{8}\,(2l + l^2)
$$
is equal to one of $-1$, $-\frac{2}{5}$, $-\frac{3}{5}$. The only solutions,
both for the eigenvalue $-1$, are
$
k_1=2,\, k_2 = 0,\, l=0 \; \mbox{and} \; k_1=0,\, k_2 = 0,\, l=2 .
$
Thus the space of infinitesimal Einstein deformations is equal to the space of infinitesimal
$\G_2$-deformations and the only two representations of $\Sp(2) \times \Sp(1)$ which could
be contained in this space are $V(2,0) \cong \mathfrak{sp}(2)$ and $V(2) \cong \mathfrak{sp}(1)$.
Next we have to determine whether these spaces admit
$H$-invariant homomorphisms to $\Lambda^3_{27}\m^*$.

If the standard representations of $\Sp(1)_u$ and $\Sp(1)_d$ are denoted by $E$ and $H$ respectively,
an arbitrary irreducible representation of $\Sp(1)_u \times \Sp(1)_d$ can
be written as $S^k E S^l H$. (In this and the next subsection we shall omit the tensor product
sign and the complexification sign). Then we have the following decompositions into irreducible
$\Sp(1)_u \times \Sp(1)_d$-representations:
$$
\begin{array}{cl}
 \mathfrak{sp}(1)_o &\cong \quad S^2 H, \qquad \mathfrak{m}' \cong  EH, \qquad \mathfrak{m} \cong  S^2 H \oplus EH,
\\[1ex]
\Lambda^3 \mathfrak{m}^* &\cong  \quad S^2 E S^2 H \oplus E S^3 H \oplus 2EH \oplus S^4 H \oplus S^2 H \oplus 2\mathbb{C},
\\[1ex]
\Lambda^3_{27} \mathfrak{m}^* &\cong  \quad S^2 E S^2 H \oplus E S^3 H \oplus EH \oplus S^4 H \oplus \mathbb{C},
\\[1ex]
V(2,0) &\cong  \quad \mathfrak{sp}(2) \cong  \quad  S^2 E \oplus EH \oplus S^2 H, \qquad V(2) \cong  \mathfrak{sp}(1) = S^2 H  .
\end{array}
$$

Since $\Lambda^3_{27} \mathfrak{m}^*$ and $\mathfrak{sp}(1)$ have no common summands, $\Hom_{\Sp(1)_u \times \Sp(1)_d}(\mathfrak{sp}(1), \Lambda^3_{27} \mathfrak{m}^*) = 0$ and therefore the $\Sp(2) \times \Sp(1)$-representation $\mathfrak{sp}(1)$ is not contained in $\Omega^3_{27}(M)$.

The only common summand of $\Lambda^3_{27} \mathfrak{m}^*$ and $\mathfrak{sp}(2)$ is $EH$, so $\Hom_{\Sp(1)_u \times \Sp(1)_d}(\mathfrak{sp}(2), \Lambda^3_{27} \mathfrak{m}^*)$ is $1$-dimensional. In order to proceed we have to find an explicit equivariant homomorphism
$A:\sp(2) \rightarrow \Lambda^3_{27} \mathfrak{m}^*$ spanning this space. Since
$*\sigma_o$ and $e^{4567}$ are the two linear independent $\Sp(1)_u \times \Sp(1)_d$-invariant
forms in $\Lambda^4\m^*$, an arbitrary embedding of $EH$ in $\Lambda^3 \mathfrak{m}^*$ is given by
$$
EH \cong \mathfrak{m}' \ni X \mapsto X \lrcorner\, (\lambda *\sigma_o + \mu e^{4567}) \in \Lambda^3 \mathfrak{m}^* .
$$
The image of this map is contained in $\Lambda^3_{27} \mathfrak{m}^*$  if and only if it is orthogonal to the $EH$ in $\Lambda^3_7 \mathfrak{m}^*$. Obviously this is equivalent to $\mu = -4\lambda$ and we can take the embedding
$$
i: EH \cong \mathfrak{m}' \to \Lambda^3_{27} \mathfrak{m}^*, \quad
EH \cong \mathfrak{m}' \ni X \mapsto X \lrcorner (*\sigma_o - 4 e^{4567}) .
$$
Hence $\Hom_{\Sp(1)_u \times \Sp(1)_d}(\mathfrak{sp}(2), \Lambda^3_{27} \mathfrak{m}^*)$ is
spanned by the equivariant homomorphism $A := i \circ p$, where $p:\mathfrak{sp}(2) \to EH$
is the orthogonal projection.

Thus $U= \sp(2)$ is the only $\Sp(2)\times \Sp(1)$-representation which remains for the
solution space of the equation $*d\varphi = -\tau_0 \varphi$, describing the infinitesimal
$\G_2$-deformations. As mentioned above, this equation is equivalent to
$
(\bar d  + \frac{5}{6}\tau_0 * ) \varphi = 0$, i.e. in the case at hand to
$(\bar d  -  2\sqrt{5}* ) \varphi = 0 $. From the results of the last part of
Section~\ref{sec-nat-red} for $V = \Lambda^3_{27} \mathfrak{m}^*$ and $U= \sp(2)$
it follows that equation (\ref{main}) with $c = -2\sqrt{5}$ must be satisfied for the
chosen $A$ and all $u \in \mathfrak{sp}(2)$. However this is not the case: take $\alpha:= e_4 \in EH
\cong \mathfrak{m}' \subset \mathfrak{sp}(2)$. Then
$$
\begin{array}{ll}
e_1 \cdot \alpha = [e_1,e_4] = -\frac{2}{\sqrt{5}}e_5, \quad & i(e_5) = 3e^{467} +  e^{137} +  e^{126} +  e^{234}, \\
e_2 \cdot \alpha = [e_2,e_4] = -\frac{2}{\sqrt{5}}e_6, \quad & i(e_6) = -3e^{457} +  e^{237} -  e^{125} -  e^{134}, \\
e_3 \cdot \alpha = [e_3,e_4] = -\frac{2}{\sqrt{5}}e_7, \quad & i(e_7) = 3e^{456} -  e^{236} -  e^{135} +  e^{124}, \\
e_4 \cdot \alpha = [e_4,e_4] = 0, \quad & i(e_4) = -3e^{567} -  e^{235} +  e^{136} -  e^{127}.
\end{array}
$$
Using these equations one easily sees that the coefficient of $e^{1234}$ in the left-hand side
of (\ref{main}) is $\frac{36}{\sqrt{5}} \neq 0$. Hence $\mathfrak{sp}(2)$ is not contained in
the space of solutions of $(\bar d - 2\sqrt{5}*) \varphi = 0$.

Since the nearly parallel $\G_2$-structure of the squashed sphere is a double covering of the
one on $\mathbb{R}P^7$ the same argument applies for the real projective space and
we obtain

\begin{Proposition}
There are no infinitesimal Einstein deformations and, in particular, no infinitesimal
$\G_2$-deformations of the nearly parallel $\G_2$-structure on the squashed sphere
$\frac{\Sp(2) \times \Sp(1)}{\Sp(1) \times \Sp(1)}$ and of the nearly parallel $\G_2$-structure
on $\mathbb{R}P^7$ inducing the second Einstein metric.
\end{Proposition}

\bigskip

%
%
\subsection{\large{The example $\frac{\SU(3) \times \SU(2)}{U(1) \times \SU(2)}$ }} ~ \\

\noindent
We denote by $\SU(2)_d$ the following embedding of $\SU(2)$ in $\SU(3) \times \SU(2)$:
$$
\SU(2)_d := \{ ( \begin{pmatrix} a & 0 \\ 0 & 1 \end{pmatrix}, a): a \in \SU(2) \}
 .
$$
The group $U(1)$ is realized as a subgroup of $\SU(3) \subset \SU(3) \times \SU(2)$ by
the embedding
$$
U(1)= \{ ( \begin{pmatrix} e^{it} & 0 & 0 \\ 0 & e^{it} & 0 \\ 0 & 0 & e^{-2it} \end{pmatrix}, 1): t \in \mathbb{R} \}  .
$$
We consider the homogeneous space $\frac{\SU(3) \times \SU(2)}{U(1) \times \SU(2)_d}$ as a normal homogeneous space taking the metric induced by $-\frac{1}{24}B$, where $B$ is the Killing form of $\mathfrak{g} = \mathfrak{su}(3) \oplus \mathfrak{su}(2)$. Then we have the reductive decomposition $\mathfrak{g} = \mathfrak{h} \oplus \mathfrak{m}$, with
$$
\mathfrak{h} = \mathfrak{u}(1) \oplus \mathfrak{su}(2)_d, \qquad \mathfrak{m} = \mathfrak{h}^\bot = \mathfrak{su}(2)_o \oplus \mathfrak{m}' .
$$
Here
$$
\mathfrak{u}(1) := \mathrm{span} \{ C \},\mbox{ where } C:=( \begin{pmatrix} i & 0 & 0 \\ 0 & i & 0 \\ 0 & 0 & -2i \end{pmatrix}, 0), \mbox{ and }
\mathfrak{su}(2)_d
:= \{ ( \begin{pmatrix} a & 0 \\ 0 & 0 \end{pmatrix}, a): a \in \mathfrak{su}(2) \}
$$
are the Lie algebras of $U(1)$ and $\SU(2)_d$ respectively and
$$
\mathfrak{su}(2)_o := \{ ( \begin{pmatrix} 2a & 0 \\ 0 & 0 \end{pmatrix}, -3a): a \in \mathfrak{su}(2) \}, \qquad
\mathfrak{m}' := \{ ( \begin{pmatrix} 0 & z \\ -\bar{z}^t & 0 \end{pmatrix}, 0): z \in \mathbb{C}^2 \}  .
$$
Let
$$
I:= \begin{pmatrix} i & 0 \\ 0 & -i \end{pmatrix} \in \mathfrak{su}(2), \quad
J:= \begin{pmatrix} 0 & -1 \\ 1 & 0 \end{pmatrix} \in \mathfrak{su}(2), \quad
K:= \begin{pmatrix} 0 & i \\ i & 0 \end{pmatrix} \in \mathfrak{su}(2) .
$$
Then we define the orientation fixing the following orthonormal frame of
$\mathfrak{m}$:
$$
e_1 := -\frac{1}{\sqrt{5}}(\begin{pmatrix} 2I & 0 \\ 0 & 0 \end{pmatrix}, -3I), \quad
e_2 := -\frac{1}{\sqrt{5}}(\begin{pmatrix} 2J & 0 \\ 0 & 0 \end{pmatrix}, -3J), \quad
e_3 := -\frac{1}{\sqrt{5}}(\begin{pmatrix} 2K & 0 \\ 0 & 0 \end{pmatrix}, -3K)  ,
$$
$$
e_4 := (\begin{pmatrix} 0 & 0 & \sqrt{2} \\ 0 & 0 & 0 \\ -\sqrt{2} & 0 & 0 \end{pmatrix}, 0), \quad
e_5 := (\begin{pmatrix} 0 & 0 & \sqrt{2}i \\ 0 & 0 & 0 \\ \sqrt{2}i & 0 & 0 \end{pmatrix}, 0),$$
$$e_6 := (\begin{pmatrix} 0 & 0 & 0 \\ 0 & 0 & \sqrt{2} \\ 0 & -\sqrt{2} & 0 \end{pmatrix}, 0), \quad
e_7 := (\begin{pmatrix} 0 & 0 & 0 \\ 0 & 0 & \sqrt{2}i \\ 0 & \sqrt{2}i & 0 \end{pmatrix}, 0)  .
$$
Then, as in the previous example we see that $\hat \T_o = \frac{2}{\sqrt{5}} \sigma_o$, where $\sigma_o$ is given by \pref{sigma}. Hence, by Lemma~\ref{hlem} we obtain a nearly parallel $\G_2$-structure on $\frac{\SU(3) \times \SU(2)}{U(1) \times \SU(2)_d}$ with $\tau_0 = -\frac{12}{\sqrt{5}}$.

Again we want to find the infinitesimal Einstein deformations of this structure. By
Corollary~\ref{Cas-eig} together with \pref{Cas-product} and \pref{Cas-SU} we are
this time looking for integers $k_1 \ge k_2 \ge k_3$ and $l_1 \ge l_2$, satisfying $-\frac{3}{2} < k_1 + k_2 + k_3 \le \frac{3}{2}$ and $-1 < l_1 + l_2 \le 1$, such that
$$
\Cas^{\SU(3) \times \SU(2)}_{V(k_1,k_2,k_3) \otimes V(l_1,l_2)} = -\frac{1}{9}(3k_1 + k_1^2 + k_2^2 - 3k_3 + k_3^2 - k_1 k_2 -k_2 k_3 - k_3 k_2) - \frac{1}{8}(2l_1 + l_1^2 - 2l_2 + l_2^2 - 2l_1 l_2)  .
$$
is equal to one of $-1$, $-\frac{2}{5}$, $-\frac{3}{5}$. The only solutions, both for
the eigenvalue $-1$, are
$$
k_1=1,\, k_2 = 0,\, k_3 = -1,\, l_1=0,\, l_2 = 0 \quad \mbox{and} \quad k_1 = 0,\, k_2 = 0,\, k_3 = 0,\, l_1 = 1,\, l_2 = -1 ,
$$
Thus the space of infinitesimal Einstein deformations is equal to the space of infinitesimal $\G_2$-deformations and the only two representations of $\SU(3) \times \SU(2)$ which could be contained in this space are $V(1,0,-1) \cong \mathfrak{su}(3)$ and $V(1,-1) \cong \mathfrak{su}(2)$. Next we have to determine the $H$-equivariant homomorphisms of these spaces into $\Lambda^3_{27}\m^*$.

If the representation of $U(1)$ with weight $k$ is denoted by $F(k)$ and the standard representation of $\SU(2)_d$ by $H$, then an arbitrary irreducible representation of $U(1) \times \SU(2)_d$ has the form $F(k) S^l H$. We have the following decompositions into irreducible $U(1) \times \SU(2)_d$-representations:
$$
\begin{array}{cl}
\mathfrak{su}(2)_o
&\cong \quad S^2 H, \qquad \mathfrak{m}' \cong F(3)H \oplus F(-3)H, \qquad \mathfrak{m} \cong S^2 H \oplus F(3)H \oplus F(-3)H,
\\[1ex]
\Lambda^3 \mathfrak{m}^*
&\cong \quad F(6) S^2 H \oplus F(-6) S^2 H \oplus F(3) S^3 H \oplus F(-3) S^3 H \oplus 2F(3) H \oplus 2F(-3) H
\\
&
\phantom{xxxxxxxxxxxxxxxxxxxxxxxxxxxxxxxxxxxxxxxxxx}
\oplus S^4 H \oplus 2S^2 H \oplus 2\mathbb{C},
\\[1ex]
\Lambda^3_{27} \mathfrak{m}^*
&\cong \quad F(6) S^2 H \oplus F(-6) S^2 H \oplus F(3) S^3 H \oplus F(-3) S^3 H \oplus F(3)H \oplus F(-3)H
\\
&
\phantom{xxxxxxxxxxxxxxxxxxxxxxxxxxxxxxxxxxxxxxxxxxx}
\oplus S^4 H \oplus S^2 H \oplus \mathbb{C},
\\[1ex]
V(1,0,-1)
&\cong \quad \mathfrak{su}(3) \cong \mathbb{C} \oplus F(3)H \oplus F(-3)H \oplus S^2 H, \qquad V(1,-1) \cong \mathfrak{su}(2) \cong S^2 H,
\end{array}
$$

The only common summand of $\Lambda^3_{27} \mathfrak{m}^*$ and $\mathfrak{su}(2)$ is $S^2 H$, so
$\Hom_{U(1) \times \SU(2)_d}(\mathfrak{su}(2), \Lambda^3_{27} \mathfrak{m}^*)$ is $1$-dimensional. Let
$$ q_2: S^2 H \cong \mathfrak{su}(2) \to \mathfrak{su}(2)_o, \quad S^2 H \cong \mathfrak{su}(2) \ni a \mapsto ( \begin{pmatrix} 2a & 0 \\ 0 & 0 \end{pmatrix}, -3a) \in \mathfrak{su}(2)_o$$
be the identification of $S^2 H$ and $\mathfrak{su}(2)_o$. The $S^2 H$ in $\Lambda^3 \mathfrak{m}^*$ coming from $\Lambda^3_7 \mathfrak{m}^*$ is given by the embedding
$$S^2 H \ni a \mapsto  q_2(a) \mapsto  q_2(a) \lrcorner *\sigma_0 \in \Lambda^3_7 \mathfrak{m}^*.$$
Since the $2$-form $\Omega := e^{45} + e^{67}$ is also $U(1) \times \SU(2)_d $-invariant, another embedding of $S^2 H$ in $\Lambda^3 \mathfrak{m}^*$ is
$$i_2: S^2 H  \to \Lambda^3 \mathfrak{m}^*, \quad  S^2 H \ni a \mapsto  q_2(a) \mapsto  q_2(a)^\flat \wedge \Omega \in \Lambda^3 \mathfrak{m}^*.$$
It is easy to see that $i_2 (S^2 H)$ is orthogonal to the $S^2 H$ in $\Lambda^3_7 \mathfrak{m}^*$, so in fact $i_2$ is the embedding of $S^2 H$ into $\Lambda^3_{27} \mathfrak{m}^*$. Therefore $\Hom_{U(1) \times \SU(2)_d}(\mathfrak{su}(2), \Lambda^3_{27} \mathfrak{m}^*)$ is spanned by $A := i_2$.

Now, as in the previous example, it remains to solve the equation $(\bar d  -2\sqrt{5}*)\varphi = 0$
by applying the results of the last part of section~\ref{sec-nat-red} for $V = \Lambda^3_{27} \mathfrak{m}^*$
and for the summands $U = \mathfrak{su}(3)$ and $U = \mathfrak{su}(2)$ found above.

If $\mathfrak{su}(2)$ is contained in the space of solutions of $(\bar d  -2\sqrt{5}*)\varphi =0$,
then equation (\ref{main}) (with $s=3$ and $c = -2\sqrt{5}$) must be satisfied for the chosen $A$ and
all $\alpha \in \mathfrak{su}(2)$. We shall show this is not the case. Take $\alpha := (0,I) \in S^2 H
\cong \mathfrak{su}(2) \subset \mathfrak{su}(3) \oplus \mathfrak{su}(2)$. Then
$$
e_2 \cdot \alpha = [e_2,\alpha] = \frac{6}{\sqrt{5}}(0,K), \quad i_2(0,K) =  q_2(0,K)^\flat \wedge
\Omega = -\sqrt{5}e^3 \wedge (e^{45} + e^{67}) = -\sqrt{5}(e^{345} + e^{367}),$$
$$e_3 \cdot \alpha = [e_3,\alpha] = -\frac{6}{\sqrt{5}}(0,J), \quad i_2(0,J) =  q_2(0,J)^\flat \wedge
\Omega = -\sqrt{5}e^2 \wedge (e^{45} + e^{67}) = -\sqrt{5}(e^{245} + e^{267}),$$
$$e_4 \cdot \alpha = [e_4,\alpha] = 0, \quad e_5 \cdot \alpha = [e_5,\alpha] = 0,$$
$$i_2(0,I) = q_2(0,I)^\flat \wedge \Omega = -\sqrt{5}e^1 \wedge (e^{45} + e^{67}) = -\sqrt{5}(e^{145} + e^{167}).$$
Using these equations one easily sees that the coefficient of $e^{2345}$ in the left-hand side of
(\ref{main}) is $22 \neq 0$. Hence $\mathfrak{su}(2)$ is not contained in the space of solutions of
$(\bar d  -2\sqrt{5}*) \varphi = 0$.

There are four common summands of $\Lambda^3_{27} \mathfrak{m}^*$ and $\mathfrak{su}(3)$: $\mathbb{C}$,
 $S^2 H$, $F(3)H$, $F(-3)H$. Since they are all different, $\Hom_{U(1) \times \SU(2)_d}(\mathfrak{su}(3),
 \Lambda^3_{27} \mathfrak{m}^*)$ is $4$-dimensional. Our next goal is to determine a basis $A_1,A_2,A_3,A_4$
 of $\Hom_{U(1) \times \SU(2)_d}(\mathfrak{su}(3), \Lambda^3_{27} \mathfrak{m}^*)$ corresponding to these
 spaces.

The $\mathbb{C}$ in $\Lambda^3 \mathfrak{m}^*$ coming from $\Lambda^3_1 \mathfrak{m}^*$ is spanned
by $\sigma_o$ and the second $\mathbb{C}$ in $\Lambda^3 \mathfrak{m}^*$ is spanned by $e^{123}$.
Thus an arbitrary $\mathbb{C}$ in $\Lambda^3 \mathfrak{m}^*$ is spanned by $\lambda \sigma_o +
\mu e^{123}$. This is orthogonal to $\sigma_o$ if and only if $\mu = -7\lambda$. Hence the
$\mathbb{C}$ in $\Lambda^3_{27} \mathfrak{m}^*$ is spanned by $\sigma_o - 7 e^{123}$. On
the other hand, $\mathbb{C}$ in $\mathfrak{su}(3)$ is $\mathfrak{u}(1)$ and is spanned by $C$. Define
$$i_1: \mathfrak{u}(1) \to \Lambda^3_{27} \mathfrak{m}^*, \quad C \mapsto \sigma_o - 7 e^{123}.$$
Then the subspace of $\Hom_{U(1) \times \SU(2)_d}(\mathfrak{su}(3), \Lambda^3_{27} \mathfrak{m}^*)$
which corresponds to $\mathbb{C}$ is spanned by $A_1 := i_1 \circ p_1$, where $p_1:\mathfrak{su}(3)
\to \mathfrak{u}(1)$ is the projection.

Let
$$j_2: S^2 H \cong \mathfrak{su}(2)  \to \mathfrak{su}(3), \quad  S^2 H \cong \mathfrak{su}(2) \ni a \mapsto  ( \begin{pmatrix} a & 0 \\ 0 & 0 \end{pmatrix}, 0) \in \mathfrak{su}(3) \subset \mathfrak{su}(3) \oplus \mathfrak{su}(2).$$
Then the subspace of $\Hom_{U(1) \times \SU(2)_d}(\mathfrak{su}(3), \Lambda^3_{27} \mathfrak{m}^*)$ corresponding to $S^2 H$ is spanned by $A_2 := i_2 \circ j_2^{-1} \circ p_2$, where $p_2:\mathfrak{su}(3) \to S^2 H$ is the projection and $i_2$ was defined earlier.

Considered as subspaces of $\mathfrak{m}' \subset \mathfrak{su}(3) \oplus \mathfrak{su}(2)$, $F(3)H$ and $F(-3)H$ are
$$F(3)H \cong \mathrm{span} \{ e_4 - ie_5, e_6 - ie_7 \}, \quad F(-3)H \cong \mathrm{span} \{ e_4 + ie_5, e_6 + ie_7 \}.$$

In the same way as for $EH$ in the case of $\frac{\Sp(2) \times \Sp(1)}{\Sp(1) \times \Sp(1)}$ we obtain that the embeddings $i_3: F(3)H \to \Lambda^3_{27} \mathfrak{m}^*$ and $i_4: F(-3)H \to \Lambda^3_{27} \mathfrak{m}^*$ are given by the restrictions on $F(3)H$ and $F(-3)H$ of the embedding
$$\mathfrak{m}' \ni X \mapsto X \lrcorner (*\sigma_o - 4 e^{4567}).$$
Then the subspaces of $\Hom_{U(1) \times \SU(2)_d}(\mathfrak{su}(3), \Lambda^3_{27} \mathfrak{m}^*)$ corresponding to $F(3)H$ and $F(-3)H$ are spanned by $A_3 := i_3 \circ p_3$ and $A_4 := i_4 \circ p_4$, where $p_3:\mathfrak{su}(3) \to F(3)H$, $p_4: \mathfrak{su}(3) \to F(-3)H$ are the projections.

Thus we have to find for which $A = c_1 A_1 + c_2 A_2 + c_3 A_3 + c_4 A_4$ equation (\ref{main}) (with $s=3$ and $c = -2\sqrt{5}$) is satisfied for all $\alpha \in \mathfrak{su}(3)$. As this equation is $U(1) \times \SU(2)_d$-invariant, this is equivalent to the requirement that the equation is satisfied for one representative of each of the four summands in $\mathfrak{su}(3)$. We take
$$
\begin{array}{ll}
\alpha_1 := C \in \mathbb{C} \subset \mathfrak{su}(3), & \alpha_2 := j_2 (I) \in S^2 H \subset \mathfrak{su}(3), \\
\alpha_3 := e_4 - ie_5 \in F(3)H \subset \mathfrak{su}(3), &
\alpha_4 := e_4 + ie_5 \in F(-3)H \subset \mathfrak{su}(3).
\end{array}
$$
Then we have
$$A(C) = c_1 i_1(C) = c_1 (\sigma_o - 7e^{123}),$$
$$A(j_2(I)) = c_2 i_2(I) = -\sqrt{5}c_2 e^1 \wedge \Omega = -\sqrt{5}c_2 (e^{145} + e^{167}),$$
$$A(j_2(J)) = c_2 i_2(J) = -\sqrt{5}c_2 e^2 \wedge \Omega = -\sqrt{5}c_2 (e^{245} + e^{267}),$$
$$A(j_2(K)) = c_2 i_2(K) = -\sqrt{5}c_2 e^3 \wedge \Omega = -\sqrt{5}c_2 (e^{345} + e^{367}),$$
$$A(e_4 - ie_5) = c_3 i_3(e_4 - ie_5) = c_3((-3e^{567} - e^{235} + e^{136} - e^{127}) - i (3e^{467} + e^{137} + e^{126} + e^{234})),$$
$$A(e_6 - ie_7) = c_3 i_3(e_6 - ie_7) = c_3((-3e^{457} + e^{237} - e^{125} - e^{134}) - i (3e^{456} - e^{236} - e^{135} + e^{124})),$$
$$A(e_4 + ie_5) = c_4 i_4(e_4 + ie_5) = c_4((-3e^{567} - e^{235} + e^{136} - e^{127}) + i (3e^{467} + e^{137} + e^{126} + e^{234})),$$
$$A(e_6 + ie_7) = c_4 i_4(e_6 + ie_7) = c_4((-3e^{457} + e^{237} - e^{125} - e^{134}) + i (3e^{456} - e^{236} - e^{135} + e^{124})).$$

Since equation (\ref{main}) is invariant with respect to $U(1) \times \SU(2)_d$, its left-hand side for $\alpha = \alpha_1$ lies in $2\mathbb{C} = \mathrm{span}\{ *\sigma_o, e^{4567} \} \subset \Lambda^4 \mathfrak{m}^*$. Hence, to determine it, it is enough to compute the coefficients of $e^{4567}$ and $e^{2367}$.

We have
$$e_1 \cdot \alpha_1 = [e_1,C] = 0, \quad e_2 \cdot \alpha_1 = [e_2,C] = 0, \quad e_3 \cdot \alpha_1 = [e_3,C] = 0,$$
$$e_4 \cdot \alpha_1 = [e_4,C] = -3e_5, \quad e_5 \cdot \alpha_1 = [e_5,C] = 3e_4,$$
$$e_6 \cdot \alpha_1 = [e_6,C] = -3e_7, \quad e_7 \cdot \alpha_1 = [e_7,C] = 3e_6.$$
Using these equations we see that the coefficients of $e^{4567}$ and $e^{2367}$ are $-18i(c_3 - c_4) + 12\sqrt{5}c_1$ and $3i(c_3 - c_4) - 2\sqrt{5}c_1$. So the whole left-hand side of (\ref{main}) for $\alpha = \alpha_1$ is
$$(-2\sqrt{5}c_1 + 3i(c_3 - c_4))(*\sigma_0 - 7e^{4567})$$
and this vanishes if and only if $-2\sqrt{5}c_1 + 3i(c_3 - c_4) = 0$.

Since equation (\ref{main}) is $U(1) \times \SU(2)_d$-invariant, its left-hand side for $\alpha \in S^2 H$ lies in $2S^2 H \subset \Lambda^4 \mathfrak{m}^*$. The two embeddings of $S^2 H$ in $\Lambda^4 \mathfrak{m}^*$ are obtained as the composition of the two embeddings of $S^2 H$ in $\Lambda^3 \mathfrak{m}^*$ with $*$. Thus the left-hand side of (\ref{main}) for $\alpha = \alpha_2$ lies in
$$\mathrm{span} \{ e^{1246} + e^{1347} + e^{1257} - e^{1356}, e^{2367} + e^{2345} \}.$$
Hence, to determine it, it is enough to compute the coefficients of $e^{1246}$ and $e^{2367}$.

We have
$$e_1 \cdot \alpha_2 = [e_1,j_2(I)] = 0, \quad e_2 \cdot \alpha_2 = [e_2,j_2(I)] = -\frac{4}{\sqrt{5}}j_2(K), \quad e_3 \cdot \alpha_2 = [e_3,j_2(I)] = \frac{4}{\sqrt{5}}j_2(J),$$
$$e_4 \cdot \alpha_2 = [e_4,j_2(I)] = -e_5, \quad e_5 \cdot \alpha_2 = [e_5,j_2(I)] = e_4,$$
$$e_6 \cdot \alpha_2 = [e_6,j_2(I)] = e_7, \quad e_7 \cdot \alpha_2 = [e_7,j_2(I)] = -e_6.$$
With these equations we see that the coefficients of $e^{1246}$ and $e^{2367}$ are $c_3 + c_4$ and $2c_2 - i(c_3 - c_4)$. So the whole left-hand side of (\ref{main}) for $\alpha = \alpha_2$ is
$$(c_3 + c_4)(e^{1246} + e^{1347} + e^{1257} - e^{1356}) + (2c_2 - i(c_3 - c_4))(e^{2367} + e^{2345})$$
and this vanishes if and only if $c_3 + c_4 = 0$ and $2c_2 - i(c_3 - c_4) = 0$.

Again the $U(1) \times \SU(2)_d$-invariance of equation (\ref{main}) implies that its left-hand side for $\alpha \in F(3)H$ lies in $2F(3)H \subset \Lambda^4 \mathfrak{m}^*$. The two embeddings of $F(3)H$ in $\Lambda^4 \mathfrak{m}^*$ are obtained as the composition of the two embeddings of $F(3)H$ in $\Lambda^3 \mathfrak{m}^*$ with $*$. Thus the left-hand side of (\ref{main}) for $\alpha = \alpha_3$ lies in
$$\mathrm{span} \{ e^{1234} - ie^{1235}, (e^{1467} - e^{2457} + e^{3456}) - i(e^{2456} + e^{3457} + e^{1567}) \}.$$
Hence, to determine it, it is enough to compute the coefficients of $e^{1234}$ and $e^{1467}$.

We have
$$e_1 \cdot \alpha_3 = [e_1,e_4 - ie_5] = -\frac{2}{\sqrt{5}}i(e_4 - ie_5), \quad e_2 \cdot \alpha_3 = [e_2,e_4 - ie_5] = -\frac{2}{\sqrt{5}}(e_6 - ie_7),$$
$$e_3 \cdot \alpha_3 = [e_3,e_4 - ie_5] = -\frac{2}{\sqrt{5}}i(e_6 - ie_7),$$
$$e_4 \cdot \alpha_3 = [e_4,e_4 - ie_5)] = -2iC - 2ij_2(I), \quad e_5 \cdot \alpha_3 = [e_5,e_4 - ie_5] = -2C - 2j_2(I),$$
$$e_6 \cdot \alpha_3 = [e_6,e_4 - ie_5] = -2j_2(J) - 2ij_2(K), \quad e_7 \cdot \alpha_3 = [e_7,e_4 - ie_5] = i(2j_2(J) + 2ij_2(K)).$$
Then, using these equations we find that the coefficients of $e^{1234}$ and $e^{1467}$ are $12ic_1 + \frac{36}{\sqrt{5}}c_3$ and $2ic_1 + 2i \sqrt{5}c_2 + \frac{16}{\sqrt{5}}c_3$. So the whole left-hand side of (\ref{main}) for $\alpha = \alpha_3$ is
$$(12ic_1 + \frac{36}{\sqrt{5}}c_3)(e^{1234} - ie^{1235})$$
$$+ (2ic_1 + 2i \sqrt{5}c_2 + \frac{16}{\sqrt{5}}c_3)((e^{1467} - e^{2457} + e^{3456}) - i(e^{2456} + e^{3457} + e^{1567}))$$
and this vanishes if and only if $12ic_1 + \frac{36}{\sqrt{5}}c_3 = 0$ and $2ic_1 + 2i \sqrt{5}c_2 + \frac{16}{\sqrt{5}}c_3 = 0$.

The computations for $\alpha = \alpha_4 \in F(-3)H$ are similar. In fact, one has to take the results for $\alpha_3$, change $c_3$ to $c_4$, preserve $c_1$ and $c_2$ and take the complex conjugate of everything else. So the whole left-hand side of (\ref{main}) for $\alpha = \alpha_4$ is
$$(-12ic_1 + \frac{36}{\sqrt{5}}c_4)(e^{1234} + ie^{1235})$$
$$+ (-2ic_1 - 2i \sqrt{5}c_2 + \frac{16}{\sqrt{5}}c_4)((e^{1467} - e^{2457} + e^{3456}) + i(e^{2456} + e^{3457} + e^{1567}))$$
and this vanishes if and only if $-12ic_1 + \frac{36}{\sqrt{5}}c_4 = 0$ and $-2ic_1 - 2i \sqrt{5}c_2 + \frac{16}{\sqrt{5}}c_4 = 0$.

Hence equation (\ref{main}) is satisfied for $A = c_1 A_1 + c_2 A_2 + c_3 A_3 + c_4 A_4$ and all $\alpha \in \mathfrak{su}(3)$ if and only if
$$-2\sqrt{5}c_1 + 3i(c_3 - c_4) = 0, \quad c_3 + c_4 = 0, \quad 2c_2 - i(c_3 - c_4) = 0,$$
$$12ic_1 + \frac{36}{\sqrt{5}}c_3 = 0, \quad 2ic_1 + 2i \sqrt{5}c_2 + \frac{16}{\sqrt{5}}c_3 = 0,$$
$$-12ic_1 + \frac{36}{\sqrt{5}}c_4 = 0, \quad -2ic_1 - 2i \sqrt{5}c_2 + \frac{16}{\sqrt{5}}c_4 = 0.$$
The solution of this linear system is $1$-dimensional:
$$c_2 = \frac{\sqrt{5}}{3}c_1, \quad c_3 = -\frac{\sqrt{5}}{3}ic_1, \quad c_4 = \frac{\sqrt{5}}{3}ic_1.$$
This means that exactly one copy of $\mathfrak{su}(3)$ is contained in the space of solutions of the equation
$(\bar d  -2\sqrt{5}*)\varphi =0$.

Thus we have proved

\begin{Proposition}
The space of infinitesimal Einstein deformations of the proper nearly parallel $\G_2$-structure on
$\frac{\SU(3) \times \SU(2)}{U(1) \times \SU(2)}$ coincides with the space of its infinitesimal
$\G_2$-deformations. This space is $8$-dimensional and is isomorphic to $\mathfrak{su}(3)$ as an
$\SU(3) \times \SU(2)$-representation.
\end{Proposition}


\begin{appendix}\label{ap1}

\end{appendix}


\end{document}